\documentclass[final,onefignum,onetabnum]{siamart220329}
\usepackage{braket,amsfonts,bm}

\usepackage{array}
\usepackage[caption=false]{subfig}
\usepackage{pgfplots}
\usepackage{cases}
\newsiamthm{claim}{Claim}
\newsiamremark{remark}{Remark}
\newsiamremark{hypothesis}{Hypothesis}
\crefname{hypothesis}{Hypothesis}{Hypotheses}

\usepackage{algorithmic}

\usepackage{graphicx,epstopdf}

\Crefname{ALC@unique}{Line}{Lines}

\usepackage{amsopn}

\usepackage{xspace}
\usepackage{bold-extra}
\usepackage[most]{tcolorbox}

\colorlet{texcscolor}{blue!50!black}
\colorlet{texemcolor}{red!70!black}
\colorlet{texpreamble}{red!70!black}
\colorlet{codebackground}{black!25!white!25}

\lstdefinestyle{siamlatex}{%
  style=tcblatex,
  texcsstyle=*\color{texcscolor},
  texcsstyle=[2]\color{texemcolor},
  keywordstyle=[2]\color{texemcolor},
  moretexcs={cref,Cref,maketitle,mathcal,text,headers,email,url},
}

\tcbset{%
  colframe=black!75!white!75,
  coltitle=white,
  colback=codebackground, 
  colbacklower=white, 
  fonttitle=\bfseries,
  arc=0pt,outer arc=0pt,
  top=1pt,bottom=1pt,left=1mm,right=1mm,middle=1mm,boxsep=1mm,
  leftrule=0.3mm,rightrule=0.3mm,toprule=0.3mm,bottomrule=0.3mm,
  listing options={style=siamlatex}
}

\newtcblisting[use counter=example]{example}[2][]{%
  title={Example~\thetcbcounter: #2},#1}

\newtcbinputlisting[use counter=example]{\examplefile}[3][]{%
  title={Example~\thetcbcounter: #2},listing file={#3},#1}

\DeclareTotalTCBox{\code}{ v O{} }
{ 
  fontupper=\ttfamily\color{black},
  nobeforeafter,
  tcbox raise base,
  colback=codebackground,colframe=white,
  top=0pt,bottom=0pt,left=0mm,right=0mm,
  leftrule=0pt,rightrule=0pt,toprule=0mm,bottomrule=0mm,
  boxsep=0.5mm,
  #2}{#1}

\patchcmd\newpage{\vfil}{}{}{}
\begin{tcbverbatimwrite}{tmp_\jobname_header.tex}
\title{The ADMM-PINNs algorithmic framework for nonsmooth PDE-constrained optimization: a deep learning approach \thanks{Submitted to the editors DATE.
\funding{The work of  Y.S was supported by the Humboldt Research Fellowship for postdoctoral
	researchers. The work of X.Y was supported by the RGC TRS project T32-707/22-N. The work of H.Y was supported by the National Natural Science Foundation of China (No. 12301399) and Natural Science Foundation of Tianjin (No. 22JCQNJC01120).}}}

\author{Yongcun Song\thanks{Chair for Dynamics, Control, Machine Learning and Numerics$-$Alexander von Humboldt-Professorship, Department of Mathematics,  Friedrich-Alexander-Universit\"at Erlangen-N\"urnberg, 91058 Erlangen, Germany
  (\email{ysong307$@$gmail.com}).}
\and Xiaoming Yuan\thanks{Department of Mathematics, The University of Hong Kong, Pok Fu Lam Road, Hong Kong,
	China
  (\email{xmyuan@hku.hk}).}
\and Hangrui Yue\thanks{School of Mathematical Sciences, Nankai University, Tianjin 300071,
	China
  (\email{yuehangrui@gmail.com}).}}

\headers{ADMM-PINNs for nonsmooth PDE-constrained optimization}{Y. Song, X. Yuan, and H. Yue}
\end{tcbverbatimwrite}
\title{The ADMM-PINNs algorithmic framework for nonsmooth PDE-constrained optimization: a deep learning approach \thanks{Submitted to the editors DATE.
\funding{The work of  Y.S was supported by the Humboldt Research Fellowship for postdoctoral
	researchers. The work of X.Y was supported by the RGC TRS project T32-707/22-N. The work of H.Y was supported by the National Natural Science Foundation of China (No. 12301399) and Natural Science Foundation of Tianjin (No. 22JCQNJC01120).}}}

\author{Yongcun Song\thanks{Chair for Dynamics, Control, Machine Learning and Numerics$-$Alexander von Humboldt-Professorship, Department of Mathematics,  Friedrich-Alexander-Universit\"at Erlangen-N\"urnberg, 91058 Erlangen, Germany
  (\email{ysong307$@$gmail.com}).}
\and Xiaoming Yuan\thanks{Department of Mathematics, The University of Hong Kong, Pok Fu Lam Road, Hong Kong,
	China
  (\email{xmyuan@hku.hk}).}
\and Hangrui Yue\thanks{School of Mathematical Sciences, Nankai University, Tianjin 300071,
	China
  (\email{yuehangrui@gmail.com}).}}

\headers{ADMM-PINNs for nonsmooth PDE-constrained optimization}{Y. Song, X. Yuan, and H. Yue}

\ifpdf
\hypersetup{ pdftitle={Guide to Using  SIAM'S \LaTeX\ Style} }
\fi

\begin{document}
\maketitle

\begin{abstract}
We study the combination of the alternating direction method of multipliers (ADMM) with physics-informed neural networks (PINNs) for a general class of nonsmooth partial differential equation (PDE)-constrained optimization problems, where additional regularization can be employed for constraints on the control or design variables. The resulting ADMM-PINNs algorithmic framework substantially enlarges the applicability of PINNs to nonsmooth PDE-constrained optimization problems. The ADMM makes it possible to untie the smooth PDE constraints and the nonsmooth regularization terms for iterations. Accordingly, at each iteration, one of the resulting subproblems is a smooth PDE-constrained optimization which can be efficiently solved by PINNs, and the other is a simple nonsmooth optimization problem which usually has a closed-form solution or can be efficiently solved by various standard optimization algorithms or pre-trained neural networks. The ADMM-PINNs algorithmic framework is mesh-free, easy to implement, scalable to different PDE settings, and does not require solving PDEs repeatedly. We validate the efficiency of the ADMM-PINNs algorithmic framework for different prototype applications, including inverse potential problems, source identification in elliptic equations, control constrained optimal control of the Burgers equation, and sparse optimal control of parabolic equations.
\end{abstract}



\section{Introduction}
Partial differential equations (PDEs) are fundamental mathematical tools for studying complex systems in physics, engineering, mechanics, chemistry, finance, and other fields. In addition to modeling and simulating such systems, it is very often to consider how to optimize or control them with certain goals. Hence, PDE-constrained optimization problems arise in various forms such as optimal control \cite{hinze2008optimization,lions1971optimal,troltzsch2010optimal}, inverse problems \cite{banks 2012,kirsch2011}, shape design \cite{lance2022}, just to mention a few.  Often, additional constraints,  such as boundedness \cite{lions1971optimal,troltzsch2010optimal}, sparsity \cite{stadler2009elliptic}, and discontinuity \cite{ChanT:2003,chavent1997regularization},  are imposed on the control variables (or design variables, depending on the application). Considerations of these additional constraints result in PDE-constrained optimization problems with nonsmooth objective functionals. Solving such a nonsmooth optimization problem is generally challenging because the nonsmoothness feature prevents straightforward applications of many canonical algorithms such as gradient descent methods, conjugate gradient methods, and quasi-Newton methods. This difficulty, together with other intrinsic difficulties such as the curse of high dimensionality and ill-conditioning of the resulting algebraic systems after mesh-based numerical discretization, makes solving such a nonsmooth PDE-constrained optimization problem extremely difficult.
Meticulous algorithm design is therefore needed. Usually, the development of an efficient algorithm requires considering the specific structure and features of the problem at hand.

\subsection{Model}
A general class of nonsmooth PDE-constrained optimization problems can be mathematically expressed as
 \begin{equation}\label{pco_model}
	\min_{y\in Y, u\in U} {J}(y,u)+ R(u)\qquad \text{s.t.}~e(y,u)=0,
 \end{equation}
where $Y$ and $U$ are Banach spaces; the operator $e: Y\times U\rightarrow Z$ with $Z$ a Banach space, and $e(y,u)=0$ represents a smooth PDE or a system of coupled smooth PDEs defined on $\bar{\Omega}=\Omega\cup\partial\Omega$ with $\Omega\subset\mathbb{R}^d  (d\geq 1)$ a bounded domain and $\partial \Omega$ its boundary. Throughout, we assume that $e(y,u)=0$ is well-posed. That is, for any given $u\in U$, it admits a unique corresponding solution $y(u)\in Y$ and $y(u)$ is continuous with respect to $u$. Above, $u$ and $y$ are called the control/design variable and the state variable, respectively. The functional $J: Y\times U\rightarrow \mathbb{R}$ consists of a data-fidelity term and a possible smooth regularization. The nonsmooth functional $R: U\rightarrow\mathbb{R}$ is employed to capture some prior information on $u$, such as sparsity, discontinuity, and lower and/or upper bounds. The goal of problem (\ref{pco_model}) is to find an unknown function $u$ such that the
objective functional $J(y,u)+R(u)$ is (locally) minimized subject to the PDE constraint $e(y,u)=0$.

Various PDE-constrained optimization problems can be covered by (\ref{pco_model}). For instance, the abstract equation $e(y,u)=0$ can correspond to different types of PDEs, such as parabolic equations \cite{glowinski1994exact,glowinski2008exact}, elliptic equations \cite{hinze2008optimization}, hyperbolic equations \cite{glowinski1995exact},  and
fractional diffusion equations \cite{biccari2021}. The unknown function $u$ can be a source term, coefficients, boundary or initial conditions for the PDE $e(y,u)=0$. The nonsmooth term $R(u)$ can be the indicator function of an admissible set, an $L^1$-regularization \cite{stadler2009elliptic}, or a total variation regularization \cite{chavent1997regularization}.
Given a target or observed function
$y_d$, typical examples covered by problem (\ref{pco_model}) include the following:

	$\bullet$ Control constrained optimal control problems \cite{hinze2008optimization,lions1971optimal,troltzsch2010optimal}: $J(y,u)=\frac{1}{2}\|y-y_d\|^2+\frac{\alpha}{2}\|u\|^2$ and $R(u)=I_{U_{ad}}(u)$ with $\alpha> 0$.  Here, $I_{U_{ad}}: U\rightarrow \mathbb{R}\cup\{\infty\}$ is the indicator function of the set $U_{ad}$ given by $I_{U_{ad}}(u)=0$ if $u\in U_{ad}$, and $I_{U_{ad}}(u)=+\infty$ otherwise. The set $U_{ad}$ can be a pointwise control set $U_{ad}=\{u\in U,| a\leq u(x)\leq b~\text{a.e.~in}~\Omega\}$ or an integral control set  $U_{ad}=\{u\in U,| a\leq \int_{\Omega}u dx\leq b\}$  with $a$ and $b$ two constants. The term $R(u)$ imposes additional boundedness constraints on the variable $u$.
	
	$\bullet$ Sparse optimal control problems \cite{biccari2022,schiindele2017,stadler2009elliptic}: $J(y,u)=\frac{1}{2}\|y-y_d\|^2+\frac{\alpha}{2}\|u\|^2$ and $R(u)=\rho\|u\|_{L^1(\Omega)}$ with $\alpha>0$ and $\rho>0$. The $L^1$-regularization $R(u)$ promotes the variable $u$ to be sparse \cite{stadler2009elliptic} and captures important applications in various fields such as optimal actuator placement \cite{stadler2009elliptic} and impulse control \cite{ciaramella2016}.
	
	$\bullet$ Discontinuous parameter or source identification \cite{ChanT:2003,ChenZ:1999}: $J(y,u)=\frac{1}{2}\|y-y_d\|^2$ and $R(u)=\gamma\int_{\Omega}|\nabla u|dx$ with $\gamma>0$ a total variation (TV) regularization. The TV-regularization $R(u)$ is capable of reserving the piecewise-constant property of $u$ and has found various applications such as image denoising and electrical impedance tomography.

\subsection{Traditional Numerical Approaches}
In the literature, there are some numerical approaches for solving some specific cases of (\ref{pco_model}). These approaches primarily focused on developing iterative schemes that can overcome the difficulty of the nonsmoothness of the underlying objective functional. For instance, the semismooth Newton (SSN) methods \cite{DeK2004,hinze2008optimization}, the inexact Uzawa method \cite{song2019}, and the alternating direction method of multipliers (ADMM) \cite{GSYY2022} have been studied within the context of control constrained optimal control problems. The SSN method \cite{stadler2009elliptic},  the primal-dual algorithm \cite{biccari2022}, and the proximal gradient method \cite{schiindele2017} have been respectively explored for solving sparse optimal control problems. The augmented Lagrangian methods \cite{ChanT:2003,ChenZ:1999} and the ADMM \cite{tyy2021} have been designed for parameter identification problems.  These approaches are based on the so-called adjoint methodology, which provides an efficient way to compute gradients. However, adjoint-based iterative schemes suffer from high computational costs due to the need for repeated solving of both the PDEs and the corresponding adjoint systems. These PDEs are typically solved by mesh-based numerical discretization schemes such as finite difference methods (FDM) or finite element methods (FEM), necessitating solutions for large-scale and ill-conditioned systems. Given that a single application of the underlying PDE solver can be very expensive, computing the PDEs and their adjoint systems repeatedly could be prohibitively expensive in practical scenarios.

\subsection{Physics-Informed Neural Networks}

In the past few years, machine learning methods such as Gaussian processes \cite{batlle2023,chen2021solving,meng2023sparse,pang2019,raissi2018numerical}, extreme learning machines \cite{calabro2021extreme,dong2021,dwivedi2020,fabiani2021,sun2019solving}, random feature methods \cite{chen2023bridging,chen2023random}, and deep learning methods \cite{cohen2023,e2018,han2018,raissi2019physics,sirignano2018dgm} have shown great power in solving various PDEs. In particular, thanks to the universal approximation property \cite{cybenko1989,hornik1991,hornik1989,leshno1993,pinkus1999} and great expressibility of deep neural networks (DNNs) \cite{lu2017,raghu2017}, deep learning has recently emerged as a powerful tool and has been used for many PDE-related problems, see e.g., \cite{beck2019,cuomo2022,e2018,han2018,karniadakis2021,khoo2021,ludeepxde2021,raissi2019physics,sirignano2018dgm,wang2021,zhu2019} and references therein. Among them, Physics-Informed Neural Networks (PINNs), which were introduced in \cite{raissi2019physics}, have been extensively studied for various PDEs. In principle, PINNs approximate the unknowns by neural networks and minimize a loss function that includes the residuals of the PDEs and initial and boundary conditions. The loss is evaluated at a set of scattered spatial-temporal points and the PDEs are solved when the loss goes to zero. In general, PINNs are mesh-free, easy to implement, and flexible to different PDEs. We refer to the review papers \cite{cuomo2022,faroughi2022,hao2022,karniadakis2021,ludeepxde2021} and references therein for more discussions.

It has been shown in \cite{mowlayi2021,raissi2019physics} that PINNs can be easily extended to PDE-constrained optimization problems by approximating the control or design variable with another neural network in addition to the one for the state variable. Then, these two neural networks are simultaneously trained by minimizing a composite loss function consisting of the objective functional and the residuals of the PDE constraint. After the training process, one can find a solution that satisfies the PDE constraint while minimizing the objective functional.  Moreover, it is worth noting that PINNs are very flexible in terms of the type of PDE, initial/boundary conditions, geometries of the domain, and objective functionals. Hence, PINNs can be applied to solve a large class of PDE-constrained optimization.  In this regard,  PINNs have been successfully applied to solve optimal control problems \cite{mowlayi2021}, inverse problems \cite{raissi2019physics}, and shape design \cite{sun2022}, just to name a few.  In addition to such direct extensions of PINNs,  the Control PINN \cite{barry2022} addresses PDE-constrained optimization by solving the first-order optimality conditions. The control/design variable, the state variable, and the adjoint variable are approximated by different DNNs. Then, by minimizing a loss function composed of the residuals of the first-order optimality conditions, a stationary point of the PDE-constrained optimization can be computed. In \cite{haoBilevel2022}, a bi-level PINNs framework is proposed for solving PDE-constrained optimization problems. In this bi-level framework, the lower-level optimization problem corresponds to the PDE constraints and the upper-level optimization problem finds a control variable for minimizing the objective functional. Minimizing the objective functional and solving the PDE constraint are thus decoupled, and the challenge of tuning the weights in the loss function of the standard PINNs is addressed.  Other PINNs-type methods for solving PDE-constrained optimization problems can be referred to \cite{basir2022,luhard2021,mitusch2021hybrid,pakravan2021solving,xu2022physics}.

\subsection{Motivation and Our Numerical Approach}

In contrast to traditional iterative methods that require discretizing the involved PDEs, PINNs-based methods are usually mesh-free, easy to implement, and flexible to different PDEs. Moreover, PINNs-based methods avoid solving adjoint systems completely by taking advantage of automatic differentiation, and could break the curse of dimensionality. Computational costs can thus be reduced substantially. Despite these advantages, it is worth noting that the above PINNs methods (except \cite{luhard2021}) only solve smooth PDE-constrained optimization, and cannot be directly applied to nonsmooth problems like (\ref{pco_model}).  A particular reason is that when the variable $u$ is approximated by a neural network $\mathcal{NN}_u$, the term $R(\mathcal{NN}_u)$ appears in the loss function. From a rigorous mathematical point of view,  the nonsmoothness of $R(\mathcal{NN}_u)$ prevents the application of commonly used neural network training technologies (e.g., back-propagation and stochastic gradient methods), which makes the neural network unable to train. The same concern also applies to other recently developed deep learning algorithms for solving PDE-constrained optimization, such as the ISMO \cite{lye2021iterative}, the operator learning methods \cite{hwang2021solving,wang2021fast}, and the amortized finite element analysis \cite{xue2020}. The ALM-based PINN in \cite{luhard2021} can treat inequality constraints but it does not apply to other nonsmooth cases like $R(u)=\|u\|_{L^1(\Omega)}$ and $R(u)=\int_{\Omega}|\nabla u|dx$ in (\ref{pco_model}). To our best knowledge, there seems no deep learning approach in the literature that can tackle the generic nonsmooth PDE-constrained optimization model (\ref{pco_model}). Hence,  it is of interest to develop new deep learning approaches for solving (\ref{pco_model}) by exploring its specific mathematical structure.

 Motivated by the aforementioned nice properties and the popularity of PINNs, we develop a PINNs-based algorithmic framework for solving (\ref{pco_model}).  Recall that it is not feasible to apply PINNs directly to solve (\ref{pco_model}). To bypass this issue, our philosophy is that the structures and properties of the model should be carefully considered in algorithmic design. One particular goal is to untie the smooth PDE constraint $e(y,u)=0$ and the nonsmooth regularization $R(u)$ so that these two inherently different terms can be treated individually.  To this end, we consider the ADMM, which is a representative operator splitting method introduced by Glowinski and Marroco in \cite{glowinski1975approximation} for nonlinear elliptic problems. The ADMM can be regarded as a splitting version of the augmented Lagrangian method (ALM) \cite{hestenes 1969,powell1969}, where the subproblem is decomposed into two parts and they are solved in the Gauss-Seidel manner. A key feature of the ADMM is that the decomposed subproblems usually are much easier than the ALM subproblems, which makes the ADMM a benchmark algorithm in various areas. In particular, the ADMM has been recently applied to solve some optimal control problems in \cite{GSY2019,GSYY2022} and inverse problems in \cite{jiao2016,tyy2021}.

In this paper, we advocate combining the ADMM with PINNs and propose the ADMM-PINNs algorithmic framework for solving  (\ref{pco_model}). The ADMM-PINNs algorithmic framework inherits all the advantages of the ADMM and PINNs. On one hand, as shown in Section \ref{se:admm-pinn}, the ADMM-PINNs algorithmic framework can treat the PDE constraint $e(y,u)=0$ and the nonsmooth regularization $R(u)$ individually, and only needs to solve two simple subproblems. One subproblem is a smooth PDE-constrained optimization problem, which can be efficiently solved by various PINNs. The other one is a simple nonsmooth optimization problem, which usually has a closed-form solution or can be efficiently solved by various well-developed optimization algorithms or pre-trained DNNs.  On the other hand,  compared with traditional numerical methods with mesh-based discretization as the PDE solvers, the ADMM-PINNs algorithmic framework does not require solving PDEs repeatedly, and it is mesh-free, easy to implement, and scalable to different types of PDEs and nonsmooth regularization.

We validate the effectiveness and flexibility of the ADMM-PINNs algorithmic framework through different case studies. In particular, we consider different cases of $e(y,u)=0$ (i.e., different PDEs) and nonsmooth regularization $R(u)$, and implement the ADMM-PINNs algorithmic framework to four prototype applications: inverse potential problem, control constrained optimal control of the Burgers equation, source identification in elliptic equations, and sparse optimal control of parabolic equations.
We also compare our results with the reference ones obtained by high-fidelity traditional numerical methods based on finite element discretization.

\subsection{Organization}
The rest of this paper is organized as follows. In Section \ref{se:admm-pinn}, we present the ADMM-PINNs algorithmic framework and discuss the solutions to the subproblems.  Then, to elaborate the implementation of the ADMM-PINNs algorithmic framework, we consider an inverse potential problem in Section \ref{se:inverse_potential}, control constrained optimal control of the Burgers equation in Section \ref{se:control_burgers},  source identification for elliptic equations in Section \ref{se:source_iden}, and sparse optimal control of parabolic equations in Section \ref{se:sparse_control}, respectively.  Some specific ADMM-PINNs algorithms are derived for solving the above-mentioned concrete applications.  The effectiveness and efficiency of the resulting ADMM-PINNs algorithms are demonstrated in each section by some preliminary numerical results. Finally, some conclusions and future perspectives are reported in Section \ref{se:conclusion}.

\section{The ADMM-PINNs Algorithmic Framework}\label{se:admm-pinn}
In this section,  we first discuss the implementation of ADMM to (\ref{pco_model}) and then take a close look at the application of PINNs to the solutions of the subproblems in each ADMM iteration. The ADMM-PINNs algorithmic framework for solving (\ref{pco_model}) is thus proposed.
\subsection{ADMM}

For any given $u\in U$, let $y(u)$  be the corresponding solution to the PDE constraint $e(y,u)=0$. Then, problem (\ref{pco_model}) can be rewritten as the following unconstrained optimization problem
\vspace{-0.5em}
 \begin{equation}\label{pco_model_e}
	\min_{u\in U} {J}(y(u),u)+ R(u).
 \end{equation}

Next, we introduce an auxiliary variable $z\in U$ satisfying $z=u$. Then, problem (\ref{pco_model_e}) is equivalent to the following linearly constrained
separable optimization problem
 \vspace{-1em}\begin{equation}\label{pco_model_e2}
	\min_{u,z \in U} J(y(u), u)+ R(z),\qquad\text{s.t.}\quad u=z.
 \end{equation}
We assume that  the canonical $L^2$-inner product $(\cdot,\cdot)$ and the $L^2$-norm$\|\cdot\|$ are well-defined in $U$. Then, an augmented Lagrangian functional associated with (\ref{pco_model_e2}) can be defined as
 \vspace{-0.5em}$$
L_{\beta}(u,z;\lambda)= {J}(y(u), u)+ R(z)-(\lambda,u-z)+\frac{\beta}{2}\|u-z\|^2,
 $$
where $\beta>0$ is a penalty parameter and $\lambda\in U$ is the Lagrange multiplier associated with $u=z$.   Then, applying the ADMM \cite{glowinski1975approximation} to (\ref{pco_model_e2}), we readily obtain the following iterative scheme
\vspace{-0.5em}\begin{subequations}\label{admm}
	\begin{numcases}
		{}u^{k+1}=\arg\min_{u\in U} L_{\beta}(u, z^k;\lambda^k)\label{admm_u},\\
		z^{k+1}=\arg\min_{z\in U} L_{\beta}(u^{k+1}, z;\lambda^k)\label{admm_z},\\
		\lambda^{k+1}=\lambda^k-\beta(u^{k+1}-z^{k+1}).
	\end{numcases}
\vspace{-0.5em}\end{subequations}

It is easy to see that the ADMM decouples the nonsmooth regularization $R(u)$ and the PDE constraint $e(y,u)=0$ at each iteration. We shall show in the rest part of this section
that the subproblem (\ref{admm_u}) can be efficiently solved by PINNs and the subproblem (\ref{admm_z}) usually has a closed-form solution or can be efficiently solved by some well-developed numerical solvers.

\subsection{Solution to the $u$-Subproblem (\ref{admm_u})}\label{se:solution_u}

We first note that the  $u$-subproblem (\ref{admm_u}) is equivalent to the following PDE-constrained optimization problem:
 \vspace{-0.5em}\begin{equation}\label{u_sub}
		\min_{y,u} \mathcal{J}^k(y,u):=J(y,u)-(\lambda^k,u-z^k)+\frac{\beta}{2}\|u-z^k\|^2\qquad \text{s.t.}~e(y,u)=0.
 \end{equation}

Next, we specify two PINNs algorithms for solving problem (\ref{u_sub}). The first one is called approximate-then-optimize (AtO) PINNs, where we first approximate problem (\ref{u_sub}) by neural networks and then apply PINNs to solve the resulting problem. The AtO-PINNs can be viewed as a generalization of the classic discretize-then-optimize approach, which first discretizes a PDE-constrained optimization by a certain numerical scheme (e.g., FDM or FEM) to obtain a discrete approximation, and then optimize the problem in the discrete setting.  More discussions on the comparison between neural network approximation and finite element discretization can be found in \cite{ludeepxde2021}. The second one is referred to as optimize-then-approximate (OtA) PINNs, which means that we derive the first-order optimality system of (\ref{u_sub}) in continuous settings and then approximate it by neural networks and solve it by PINNs.   The OtA-PINNs is a generalization of the classic optimize-then-discretize approach, which discretizes the continuous optimality system by a certain numerical scheme.

\subsubsection{An Approximate-then-Optimize PINNs Algorithm}
In this subsection, by generalizing the classic discretize-then-optimize approach, we specify an  AtO-PINNs algorithm to solve (\ref{u_sub}).

To approximate problem (\ref{u_sub}), we construct two neural networks $ \hat{y}(x; \bm{\theta}_y)$ parameterized by $\bm{\theta}_y$ and $\hat{u}(x;\bm{\theta}_u)$ parameterized by $\bm{\theta}_u$ to approximate $y$ and $u$, respectively.
We can show that subproblem (\ref{u_sub}) can be approximated by the following optimization problem:
  \begin{equation}\label{pco_model_admm}
	\small{\left\{
	\begin{aligned}
		&\min_{\bm{\theta}_y,\bm{\theta}_u}  J(\hat{y}(x; \bm{\theta}_y),\hat{u}(x; \bm{\theta}_u))-(\lambda^k(x),\hat{u}(x; \bm{\theta}_u)-z^k(x))+\frac{\beta}{2}\|\hat{u}(x; \bm{\theta}_u)-z^k(x)\|^2 \\
		&\qquad \text{s.t.}~e(\hat{y}(x; \bm{\theta}_y),\hat{u}(x; \bm{\theta}_u))=0.
	\end{aligned}
	\right.}
  \end{equation}

We then apply the vanilla PINNs \cite{raissi2019physics} to solve problem (\ref{pco_model_admm}) and the resulting algorithm is listed in Algorithm \ref{alg:pinn_d}.
 \begin{algorithm}[htpb]
	\caption{AtO-PINNs  for (\ref{u_sub})}\label{alg:pinn_d}
	\begin{algorithmic}[1]
		\STATE{\textbf{Input:}} Initial parameters $ \bm{\theta}_{y}^0,\bm{\theta}_{u}^0$ and weights $w_e>0, w_o>0$
		\STATE  Choose training set $\mathcal{T}=\mathcal{T}_i\cup \mathcal{T}_b$ with $ \mathcal{T}_i\subset\Omega$ and $\mathcal{T}_b\subset\partial\Omega$.
		\STATE Train the neural networks to find $\bm{\theta}_y^{k+1}$ and $\bm{\theta}_u^{k+1}$  by minimizing the following total
		loss function
		 \begin{equation}\label{loss_pco_admm}
			\mathcal{L}_{total}^k(\bm{\theta}_y,\bm{\theta}_u)=w_o\mathcal{J}^k(\bm{\theta}_y,\bm{\theta}_u)+w_e\mathcal{L}_{PDE}(\bm{\theta}_y,\bm{\theta}_u),
		 \end{equation}
	where
	{\footnotesize $$\mathcal{J}^k(\bm{\theta}_y,\bm{\theta}_u):=  \frac{1}{|\mathcal{T}|}\sum_{x\in \mathcal{T}}\left(J(\hat{y}(x; \bm{\theta}_y),\hat{u}(x; \bm{\theta}_u))-(\lambda^k(x),\hat{u}(x; \bm{\theta}_u)-z^k(x))+\frac{\beta}{2}\|\hat{u}(x; \bm{\theta}_u)-z^k(x)\|^2\right), $$}
	and
 $$ \mathcal{L}_{PDE}(\bm{\theta}_y,\bm{\theta}_u)=\frac{1}{|\mathcal{T}|}\sum_{x\in \mathcal{T}}\|e(\hat{y}(x; \bm{\theta}_y),\hat{u}(x;\bm{\theta}_u))\|^2. $$
		\STATE{\textbf{Output:}} Parameters $(\bm{\theta}_{y}^{k+1},\bm{\theta}_{u}^{k+1})$ and   solutions $\hat{y}(x;\bm{\theta}^{k+1}_{y})$ and $\hat{u}(x;\bm{\theta}^{k+1}_{u})$.
	\end{algorithmic}
\end{algorithm}

\begin{remark}
	To expose our main ideas clearly, the vanilla PINNs algorithm is employed in Algorithm \ref{alg:pinn_d} for solving (\ref{u_sub}). It is worth noting that the Deep Galerkin method \cite{sirignano2018dgm}, the Deep Ritz method \cite{e2018}, and some other PINNs-based algorithms (e.g., VPINN \cite{kharazmiVPINN}, gPINN \cite{yugPINN}, and CAN-PINN \cite{chiu_canPINN}) can also be applied in a similar way.
\end{remark}

        In the context of PINNs, the training sets $\mathcal{T}_i$ and $\mathcal{T}_b$ are usually referred to as the sets of residual points.
        For the choice of the residual points, there are generally three possible strategies: 1) specify the residual points at the beginning of training, which could be grid points on a lattice or random points, and never change them during the training process; 2) select different residual points randomly in each optimization iteration of the training process; 3) improve the location of the residual points adaptively during the training process, see e.g., \cite{ludeepxde2021}.
	  Moreover, we remark that for time-dependent problems, the time variable $t$ can be treated as an additional space coordinate and $\Omega$ denotes the spatial-temporal domain. Then, the initial and boundary conditions can be treated  similarly. Note that in Algorithm \ref{alg:pinn_d}, the calculation of $\mathcal{L}_{total}^k(\bm{\theta}_y,\bm{\theta}_u)$ usually involves derivatives, such as the partial derivatives $\frac{\partial \hat{y}}{\partial x}$ and $\frac{\partial^2 \hat{y}}{\partial x^2}$, which are handled via automatic differentiation.

Starting from some initialized parameters $\bm{\theta}_y^0$ and $\bm{\theta}_u^0$, we use stochastic optimization algorithms to find $\bm{\theta}^{k+1}_y$ and $\bm{\theta}_u^{k+1}$ that minimize (\ref{loss_pco_admm}). For instance, at each inner
iteration indexed by $j$, the parameters from both neural networks can concurrently be updated by
  $$
\begin{aligned}
	\bm{\theta}_y\leftarrow\bm{\theta}_y-\eta_j\frac{\partial \mathcal{L}^k_{total}}{\partial \bm{\theta}_y}(\bm{\theta}_y,\bm{\theta}_u),~\text{and}~
	\bm{\theta}_u\leftarrow\bm{\theta}_u-\eta_j\frac{\partial \mathcal{L}^k_{total}}{\partial \bm{\theta}_u}(\bm{\theta}_y,\bm{\theta}_u),
\end{aligned}
  $$
where $\eta_j>0$ is an appropriate learning rate. All gradients (i.e. $\frac{\partial \mathcal{L}^k_{total}}{\partial \bm{\theta}_y}$ and $\frac{\partial \mathcal{L}^k_{total}}{\partial \bm{\theta}_u}$) are computed using automatic differentiation. At the end of the training process,
the trained neural networks $\hat{y}(x;\bm{\theta}_y^{k+1})$ and  $\hat{u}(x;\bm{\theta}_u^{k+1})$ approximately solve problem (\ref{pco_model_admm}).

\subsubsection{An Optimize-then-Approximate PINNs Algorithm}
In this subsection, we derive an OtA-PINNs algorithm to solve the $u$-subproblem (\ref{u_sub}).  A similar idea can also be found in \cite{barry2022}.

Following some standard arguments as those in \cite{hinze2008optimization}, one can show that, if  $({y}, u)$ is a solution of (\ref{u_sub}), the optimality system of (\ref{u_sub}) reads as:
  \begin{equation}\label{oc_u}
\left\{
\begin{aligned}
e_u ({y},u)^* p +\mathcal{J}^k_u({y},u)=0,\\
e({y},u)=0,\\
e_y({y},u)^*p + \mathcal{J}^k_y({y},u)=0,
\end{aligned}
\right.
  \end{equation}
where $p$ is the corresponding adjoint variable, $\mathcal{J}_u^k$ and $\mathcal{J}_y^k$ are the first-order partial derivatives of $\mathcal{J}^k$ with respect to $u$ and $y$, respectively . Then,  one can compute a solution $(y,u)$ of problem (\ref{u_sub}) by solving the optimality system (\ref{oc_u}). For this purpose,  we construct three neural networks $ \hat{y}(x; \bm{\theta}_y)$ parameterized by $\bm{\theta}_y$, $\hat{u}(x;\bm{\theta}_u)$ parameterized by $\bm{\theta}_u$, and $ \hat{p}(x; \bm{\theta}_p)$ parameterized by $\bm{\theta}_p$, to approximate $y$, $u$ and $p$, respectively. Then, we approximate the optimality system (\ref{oc_u})  by
  \begin{equation}\label{oc_u_app}
	\left\{
	\begin{aligned}
		e_u (\hat{y}(x; \bm{\theta}_y),\hat{u}(x; \bm{\theta}_u))^* \hat{p}(x; \bm{\theta}_p) +\mathcal{J}^k_u(\hat{y}(x; \bm{\theta}_y),\hat{u}(x;\bm{\theta}_u))=0,\\
		e(\hat{y}(x; \bm{\theta}_y),\hat{u}(x;\bm{\theta}_u))=0,\\
		e_y(\hat{y}(x; \bm{\theta}_y),\hat{u}(x;\bm{\theta}_u))^* \hat{p}(x; \bm{\theta}_p) + \mathcal{J}^k_y(\hat{y}(x; \bm{\theta}_y),\hat{u}(x;\bm{\theta}_u))=0.
	\end{aligned}
	\right.
 \end{equation}

A PINNs algorithm for solving (\ref{oc_u_app}) is presented in Algorithm \ref{alg:pinn_oc}.
 \begin{algorithm}[htpb]
	\caption{OtA-PINNs for (\ref{u_sub})}\label{alg:pinn_oc}
	\begin{algorithmic}[1]
		\STATE{\textbf{Input:}} Initial parameters $ \bm{\theta}_{y}^0,\bm{\theta}_{u}^0, \bm{\theta}_{p}^0$ and weights $w_u>0, w_y>0, w_p>0$.
		\STATE Choose training sets $\mathcal{T}=\mathcal{T}_i\cup \mathcal{T}_b$ with $ \mathcal{T}_i\subset\Omega$ and $\mathcal{T}_b\subset\partial\Omega$.
		\STATE Specify a loss function by summing the weighted $L_2$ norm
residuals of  (\ref{oc_u_app}).
		{\small \begin{equation}\label{loss_oc}
			\begin{aligned}
				\mathcal{L}_{OS}(\bm{\theta}_y,&\bm{\theta}_u,\bm{\theta}_p)=\frac{w_y}{|\mathcal{T}|}\sum_{x\in \mathcal{T}}\|e(\hat{y}(x; \bm{\theta}_y),\hat{u}(x;\bm{\theta}_u))\|^2\\
				&+\frac{w_u}{|\mathcal{T}|}\sum_{x\in \mathcal{T}}\|e_u (\hat{y}(x; \bm{\theta}_y),\hat{u}(x; \bm{\theta}_u))^* \hat{p}(x; \bm{\theta}_p) +\mathcal{J}^k_u(\hat{y}(x; \bm{\theta}_y),\hat{u}(x;\bm{\theta}_u))\|^2\\
				&+\frac{w_p}{|\mathcal{T}|}\sum_{x\in \mathcal{T}}\|e_y(\hat{y}(x; \bm{\theta}_y),\hat{u}(x;\bm{\theta}_u))^* \hat{p}(x; \bm{\theta}_p) + \mathcal{J}^k_y(\hat{y}(x; \bm{\theta}_y),\hat{u}(x;\bm{\theta}_u))\|^2.
			\end{aligned}
		 \end{equation}}
		\STATE Train the neural networks to find $\bm{\theta}_y^{k+1}$, $\bm{\theta}_u^{k+1}$, and $\bm{\theta}_p^{k+1}$  by minimizing  (\ref{loss_oc}).
		\STATE{\textbf{Output:}} Parameters $(\bm{\theta}_{y}^{k+1},\bm{\theta}_{u}^{k+1})$ and solutions $\hat{y}(x;\bm{\theta}^{k+1}_{y})$ and $\hat{u}(x;\bm{\theta}^{k+1}_{u})$.
	\end{algorithmic}
\end{algorithm}

\begin{remark}
	Instead of the full optimality system (\ref{oc_u}), we can also use a reduced optimality system obtained by eliminating $u, y$ or $p$ in (\ref{oc_u}) to implement Algorithm \ref{alg:pinn_oc}. More details are presented in Sections  \ref{se:control_burgers} and \ref{se:source_iden}.
\end{remark}

\subsection{Solution to the $z$-Subproblem (\ref{admm_z})}\label{se:solution_z}
  It can be shown that subproblem (\ref{admm_z}) is a simple optimization problem corresponding to the nonsmooth term $R(z)$:
 \begin{equation}\label{sub_z}
z^{k+1}=\arg\min_{z}  R(z)+\frac{\beta}{2}\|z-(u^{k+1}-\frac{\lambda^k}{\beta})\|^2,
 \end{equation}
where $u^{k+1}:=\hat{u}(x;\bm{\theta}_u^{k+1})$ is the solution of the $u$-subproblem (\ref{admm_u}). Problem (\ref{sub_z}) usually has a closed-form solution or can be efficiently solved by some well-developed numerical solvers. For instance,

$\bullet$ if $R(z)=I_{U_{ad}}(z)$ is the indicator function of the set $U_{ad}=\{z\in U,| a\leq z(x)\leq b~\text{a.e. in}~\Omega\}$ with $a$ and $b$ two constants, then $z^{k+1}(x)=\mathbb{P}_{U_{ad}}\left(u^{k+1}-\frac{\lambda^k}{\beta}\right)(x)$,  a.e. in $\Omega$, where
$\mathbb{P}_{U_{ad}}(z):=\max\{\min\{z,b\},a\}$, $\forall z\in L^2(\Omega)$, is the projection onto $U_{ad}$.
	
$\bullet$ if $R(z)=\rho\|z\|_{L^1(\Omega)}$ is an $L^1$-regularization term, then {\small$z^{k+1}(x)=\mathbb{S}_{\frac{\rho}{\beta}}\left(u^{k+1}-\frac{\lambda^k}{\beta}\right)(x),$} a.e. in $ \Omega$, with $\mathbb{S}$ the Shrinkage operator defined by
	 \begin{equation}\label{def:shrinkage}
		\mathbb S_\zeta(v)(x) = \text{sgn}(v(x)) (|v(x)|-\zeta)_{+}~\text{for any}~\zeta>0,
	 \end{equation}
where $\text{sgn}$ is the sign function and $(\cdot)_{+}$ denotes the positive part.

$\bullet$ if $R(z)=\gamma \int_{\Omega}|\nabla z|dx$ is a TV-regularization term defined by  \begin{equation}\label{eq:TV}
    	\int_{\Omega}|\nabla z|dx:=\sup\Big\{\int_{\Omega}z~\mathrm{div} \varphi~\mathrm{d}x: \varphi\in C_c^{1}(\Omega;\mathbb{R}^d),\|\varphi\|_{\infty}\le1\Big\},
     \end{equation}
    where $\|\varphi\|_{\infty}=\sup_{x\in\Omega}(\sum_{i=1}^d|\varphi_i(x)|^2)^{1/2}$, ``div'' denotes the divergence operator, and $C_c^{1}(\Omega;\mathbb{R}^d)$ is the set of once
    continuously differentiable $\mathbb{R}^d$-valued functions with compact support in $\Omega$, see, e.g., \cite{Ziemer:1989} for more details. In this case, the resulting $z$-subproblem (\ref{sub_z}) has no closed-form solution, but iterative numerical algorithms for solving such a problem have been extensively studied in the literature, such as the split Bregman method \cite{css2009}, the primal-dual method \cite{chan1999}, and Chambolle’s dual method \cite{chambolle2004}. Moreover, when $\Omega\subset\mathbb{R}^2$ is a rectangle,  the $z$-subproblem (\ref{sub_z}) can be viewed as an image denoising model, which, as shown in \cite{tyy2021}, can be efficiently solved by some pre-trained deep convolutional neural networks (CNNs), see Section \ref{se:source_iden} for the details.  Here, we present an ADMM algorithm, which is closely related to the split Bregman method \cite{css2009},  to solve problem (\ref{sub_z}).  For this purpose, we first reformulate (\ref{sub_z}) as the following constrained optimization problem:
     \begin{equation}\label{sub_z_eq}
    	\min \gamma\int_{\Omega}|w|dx+\frac{\beta}{2}\|z-(u^{k+1}-\frac{\lambda^k}{\beta})\|^2,~\text{s.t.}~\nabla z=w.
     \end{equation}
    An augmented Lagrangian functional associated with problem (\ref{sub_z_eq}) is defined as
     $$
    L_{\zeta}(z,w;\mu)=\gamma\int_{\Omega}|w|dx+\frac{\beta}{2}\|z-(u^{k+1}-\frac{\lambda^k}{\beta})\|^2-(\mu, \nabla z-w)+\frac{\zeta}{2}\|\nabla z-w\|^2,
     $$
    where $\zeta>0$ is a penalty parameter and $\mu$ is the Lagrange multiplier associated with $\nabla z=w$.
    Then the iterative scheme of ADMM for solving problem (\ref{sub_z_irc_e}) reads as
    \begin{subequations}
    	{ \begin{numcases}
    		{}z_{j+1}=\arg\min \frac{\beta}{2}\|z-(u^{k+1}-\frac{\lambda^k}{\beta})\|^2+\frac{\zeta}{2}\|\nabla z-(w_j+\frac{\mu_j}{\zeta})\|^2\label{sub_z_admm1},\\
    		w_{j+1}=\arg\min\gamma\int_{\Omega}|w|dx+\frac{\zeta}{2}\|w-(\nabla z_{j+1}-\frac{\mu_j}{\zeta})\|^2\label{sub_z_admm2},\\
    		\mu_{j+1}=\mu_j-\zeta(\nabla z_{j+1}-w_{j+1})\label{sub_z_admm3}.
    	\end{numcases}}
    \end{subequations}
    The subproblems (\ref{sub_z_admm1}) and (\ref{sub_z_admm2}) can be efficiently solved. Precisely, the first-order optimality condition of (\ref{sub_z_admm1}) is a linear elliptic equation, which can be cheaply solved by multigrid preconditioned conjugate gradient methods. For (\ref{sub_z_admm2}), it has a closed form solution, which is given by
     $$w_{j+1}=\mathbb{S}_{\frac{\gamma}{\zeta}}\left(\nabla z_{j+1}-\frac{\mu_j}{\zeta}\right), $$
    where $\mathbb{S}$ is the Shrinkage operator defined by (\ref{def:shrinkage}).


\subsection{The ADMM-PINNs Algorithmic Framework for Solving (\ref{pco_model})}
Based on the discussions in Sections \ref{se:solution_u} and \ref{se:solution_z}, we propose the ADMM-PINNs algorithmic framework for problem (\ref{pco_model}) and list it in Algorithm \ref{alg:admm_pinn}.
 \begin{algorithm}[htpb]
	\caption{The ADMM-PINNs algorithmic framework for (\ref{pco_model})}\label{alg:admm_pinn}
	\begin{algorithmic}[1]
		\STATE{\textbf{Input:}} $\beta>0$, $z^0, \lambda^0, \bm{\theta}_{y}^0,\bm{\theta}_{u}^0$.
		\FOR{$k\geq 1$}
		\STATE Update $(\bm{\theta}_{y}^{k+1},\bm{\theta}_{u}^{k+1})$ and hence $ \hat{y}(x; \bm{\theta}_y^{k+1})$ and $ \hat{u}(x; \bm{\theta}_u^{k+1})$ by \textbf{Algorithm \ref{alg:pinn_d}} or by \textbf{Algorithm \ref{alg:pinn_oc}}.
	\STATE Update $z^{k+1}(x)$ by solving problem (\ref{sub_z}).
	\STATE $\lambda^{k+1}(x)=\lambda^k(x)-\beta(\hat{u}(x;\bm{\theta}_u^{k+1})-z^{k+1}(x)).$
		\ENDFOR
		\STATE{\textbf{Output:}} Parameters $(\bm{\theta}_{y}^*,\bm{\theta}_{u}^*)$ and approximate solutions $\hat{y}(x;\bm{\theta}^*_{y})$ and $\hat{u}(x;\bm{\theta}^*_{u})$.
	\end{algorithmic}
\end{algorithm}

As shown in the following sections, the  ADMM-PINNs algorithmic framework is feasible for a wide range of nonsmooth PDE-constrained optimization modeled by (\ref{pco_model}). It inherits all the advantages of the ADMM and PINNs, and has the following favorable properties. From the ADMM perspective, the nonsmooth regularization $R(u)$ is untied from the main PDE-constrained optimization problem and thus can be treated individually. Consequently, a nonsmooth PDE-constrained optimization problem is decoupled into two easier subproblems that can be efficiently solved by existing numerical methods. From the PINNs perspective, the algorithmic framework avoids solving PDEs repeatedly.  Meanwhile, it is very flexible in terms of the type of PDEs, boundary conditions, geometries, and objective functionals, and can be easily implemented to a very large class of problems without much effort on algorithmic design and coding.

 With the above-mentioned nice properties, it is promising to implement Algorithm \ref{alg:admm_pinn} to various concrete applications and propose some specific ADMM-PINNs algorithms.  To this end, we consider four classic and important inverse and optimal control problems in the following sections. Recall that we can solve the resulting $u$-subproblem (\ref{admm_u}) by Algorithm \ref{alg:pinn_d} or Algorithm \ref{alg:pinn_oc}. As a result,  two different ADMM-PINNs algorithms can be specified for each problem, and we call them ADMM-AtO-PINNs and ADMM-OtA-PINNs, respectively.


In each section, numerical results are also presented to validate the effectiveness of the specified ADMM-PINNs algorithms. All  codes used in this manuscript are written in Python and PyTorch, and are publicly available on GitHub at: \url{https://github.com/yuehangrui/ADMM_PINNs}
\section{Inverse Potential Problem}\label{se:inverse_potential}
In this section, we apply the ADMM-PINNs algorithmic framework to recover the potential term for an elliptic equation.  Such a problem plays a crucial role in different contexts of applied sciences, e.g.,  damping design, identifying heat
radiative coefficient, aquifer analysis, and optical tomography. We refer the reader to \cite{Tarantola2005} for a survey on this subject.

\subsection{Problem Statement}Let $\Omega$ be a bounded domain in $\mathbb{R}^d$ ($d\geq 1$) with a piecewise smooth boundary $\partial \Omega$. Consider the following elliptic equation
 \begin{equation}\label{eq:irc}
	\begin{aligned}
		-\nu\Delta y+uy=f \text{~in~} \Omega, y=0  \text{~on~} \partial \Omega,
	\end{aligned}
 \end{equation}
where $f\in H^{-1}(\Omega)$ is given.   Suppose that
some measurements of $y$, denoted by $y^{\delta}$, are available subject to
some noise with the noisy level $\delta > 0$. The diffusion coefficient $\nu$ in (\ref{eq:irc}) is assumed to be known and constant, but the potential $u$ is unknown and needs to be identified from the measurements $y^{\delta}$.
Here, we are interested in identifying a discontinuous potential $u$, which can be modeled by the following  TV-regularized PDE-constrained optimization problem:
 \begin{equation}\label{model:irc}
	\min_{y\in Y,u\in U} \frac{1}{2}\int_\Omega | y- y^{\delta}|^2 dx+\gamma \int_{\Omega}|\nabla u|dx, ~\text{s.t.}  -\nu\Delta y+u y=f \text{~in~} \Omega, ~y=0  \text{~on~} \partial \Omega, 
 \end{equation}
where $\gamma>0$, $Y=L^2(\Omega)$, and $U=L^2(\Omega)\cap BV(\Omega)$. It follows from \cite{chavent1997regularization} that $U$ is a Banach space endowed with the norm $\|u\|_{U}=\|u\|_{L^2(\Omega)}+\|u\|_{BV(\Omega)}$ and $U=BV(\Omega)$ when $d=1$ or $2$. The TV regularization $\int_{\Omega}|\nabla u|dx$ defined by (\ref{eq:TV})  has been widely used in inverse problems because it is capable of reserving the piecewise-constant property of $u$, see e.g., \cite{ChanT:2003,chavent1997regularization,ChenZ:1999}.

\subsection{ADMM-PINNs for (\ref{model:irc})}
Let $y(u)$ be the solution of equation (\ref{eq:irc}) corresponding to $u$. By introducing an auxiliary variable $z$ such that $u=z$, we can reformulate problem (\ref{model:irc}) as
 \begin{equation}\label{model:irc_e}
	\begin{aligned}
		\min_{u,z}&~\frac{1}{2}\int_\Omega | y(u)- y^{\delta}|^2\mathrm{d}x+\gamma\int_{\Omega}|\nabla z|dx,\quad
		\text{s.t.}~&~u=z.
	\end{aligned}
 \end{equation}
Since $U= L^2(\Omega)\cap BV(\Omega)$, the $L^2$-inner product and the $L^2$-norm are well-defined in $U$. Hence, the augmented Lagrangian functional associated with (\ref{model:irc_e}) can be defined as
 \begin{equation}\label{L:irc_e}
L_{\beta}(u,z;\lambda)= \frac{1}{2}\int_\Omega | y(u)- y^{\delta}|^2\mathrm{d}x+\gamma\int_{\Omega}|\nabla z|dx-(\lambda,u-z)+\frac{\beta}{2}\|u-z\|^2,
 \end{equation}
 where $\beta>0$ is a penalty parameter. Then the ADMM algorithm for (\ref{model:irc}) follows from \eqref{admm} with $L_{\beta}(u,z;\lambda)$ defined by \eqref{L:irc_e}.
We first note that the $z$-subproblem (\ref{admm_z}) is a TV-regularized optimization problem:
 \begin{equation}\label{sub_z_irc_e}
	z^{k+1}=\arg\min_{z} \gamma\int_{\Omega}|\nabla z|dx-(\lambda^k,u^{k+1}-z)+\frac{\beta}{2}\|u^{k+1}-z\|^2,
 \end{equation}
which can be solved by (\ref{sub_z_admm1})-(\ref{sub_z_admm3}).
The $u$-subproblem (\ref{admm_u}) is equivalent to the following PDE-constrained optimization problem:
 $$
\min \mathcal{J}^k_{IRC}(y,u):=\frac{1}{2}\int_\Omega | y- y^{\delta}|^2\mathrm{d}x-(\lambda^k, u-z^k)+\frac{\beta}{2}\|u-z^k\|^2,
 $$
subject to the equation (\ref{eq:irc}).  Its first-order optimality system reads as:
 \begin{equation}\label{oc_u_reaction}
	\left\{
	\begin{aligned}
		&p +\beta u-\lambda^k-\beta z^k=0,\\
		& -\nu\Delta y+u y=f \text{~in~} \Omega,~ y=0 \text{~on~} \partial \Omega,\\
		&-\nu \Delta p+up= y-y^{\delta} \text{~in~} \Omega,~ p=0 \text{~on~} \partial \Omega,
	\end{aligned}
	\right.
 \end{equation}
where $u$ is a solution to the $u$-subproblem (\ref{admm_u}), $y$ is the corresponding solution of  (\ref{eq:irc}), and $p$ is the adjoint variable.

 As discussed in Section \ref{se:solution_u},  an ADMM-AtO-PINNs and an ADMM-OtA-PINNs can be specified from Algorithm \ref{alg:admm_pinn} for (\ref{model:irc}).

\begin{remark}
 In our implementation of the ADMM-AtO-PINNs and  the ADMM-OtA-PINNs, the boundary conditions are enforced as soft constraints in the loss functions, which can be used for complex
	domains and any type of boundary conditions. On the other hand, it is possible to
	enforce the boundary conditions as hard constraints for some specific cases. For instance, when $\Omega  = (0, 1)$ and the boundary condition is $y(0) = y(1) = 0$, we can approximate $y$ by $\hat{y}(x;\bm{\theta}_y)=x(x-1)\mathcal{NN}_y$ with $\mathcal{NN}_y$ a neural network so that the boundary condition can be satisfied automatically; see \cite{ludeepxde2021,luhard2021} for more details.
\end{remark}

\subsection{Numerical Results for (\ref{model:irc})}
In this subsection, we test an inverse potential problem to verify the effectiveness and efficiency of the ADMM-AtO-PINNs and the ADMM-OtA-PINNs.

\noindent\textbf{Example 1.}
We set $\Omega=(0,1)$, $\nu=5\times 10^{-3}$,  {$\gamma=8\times10^{-3}$}, and $f=\sin(2\pi x), x\in \Omega$ in (\ref{model:irc}).
We further take the true potential term  $u$ as the following piece-constant function:
 $$
u_{true}=\left \{
\begin{aligned}
&1, &x\in [0.25, 0.75],\\
&0.2,& x\in (0, 0.25)\cup (0.75, 1).
\end{aligned}
\right.
 $$
The observed data $y^\delta=y_{true}+\delta \| y_{true}\| \text{rand} (x)$, where $\delta=0.05$, $y_{true}$ is the solution of \eqref{eq:irc} corresponding to $u_{true}$, and $\text{rand} (x)$ is a standard normal distributed random function.  In our numerical experiments, $y_{true}$ is the solution of \eqref{eq:irc} obtained from $u_{true}$ by FEM. To implement the FEM used in this example, the domain $\Omega=(0,1)$ is uniformly partitioned with mesh size $h=1/100$.

For the implementation of PINNs,  we approximate $u$ by a fully-connected feed-forward neural network $\hat{u}(x,\bm{\theta}_u)$. Moreover, $y$ and $p$ are respectively approximated by $\hat{y}(x;\bm{\theta}_y)=x(x-1)\mathcal{NN}(x;\bm{\theta}_y)$ and $\hat{p}(x;\bm{\theta}_p)=x(x-1)\mathcal{NN}(x;\bm{\theta}_p)$, with $\mathcal{NN}(x;\bm{\theta}_y)$ and $\mathcal{NN}(x;\bm{\theta}_p)$ being fully-connected feed-forward neural networks, so that the boundary conditions $y(0)=y(1)=0$ and $p(0)=p(1)=0$ can be satisfied automatically. Hence, we set $w_b=0$. All other weights in the loss function are set to be 1.  All the neural networks consist of 4 hidden layers with 50 neurons per hidden layer and hyperbolic tangent activation functions.  To train the neural networks, we take the training set $\mathcal{T}_i$ as the grid nodes of the lattice used in the FEM.  We first employ 20000 iterations of the Adam \cite{kingma2015} with learning rate $\eta=10^{-4}$ to obtain good guesses of the parameters and then switch to the L-BFGS \cite{nocedal1980} with a strong Wolfe line search strategy for 1000 iterations to achieve higher accuracy.  The neural networks are initialized using the default initializer of Pytorch.

For solving \eqref{sub_z_irc_e}, we discretize it by FEM and then apply the discrete version of the ADMM (\ref{sub_z_admm1})-(\ref{sub_z_admm3}). For simplicity, we assume the function $z$ in \eqref{sub_z_admm1} satisfies the periodic-Dirichlet boundary condition $z(0)=z(1)$. As a result, the discrete problem of \eqref{sub_z_admm1} can be directly solved via discrete Fast Fourier transform. We take $\zeta=0.5$ and execute 80 iterations of (\ref{sub_z_admm1})-(\ref{sub_z_admm3}) to obtain an approximate solution of \eqref{sub_z_irc_e}.

Moreover, to validate the effectiveness of the ADMM-AtO-PINNs and the ADMM-OtA-PINNs, we compare our results with the ones obtained by high-fidelity traditional numerical methods based on FEM. To be concrete, the $u$-subproblem (\ref{admm_u}) is first discretized by FEM and then is solved by a Gauss-Newton method with backtracking line search. An ADMM-Gauss-Newton method is thus specified for solving (\ref{model:irc}).  

\begin{figure}[htpb]
	\caption{\small Numerical results for Example 1}\label{fig:reaction_id}
	\centering
	\subfloat[Recovered coefficients and errors by different algorithms]{\includegraphics[width=\textwidth]{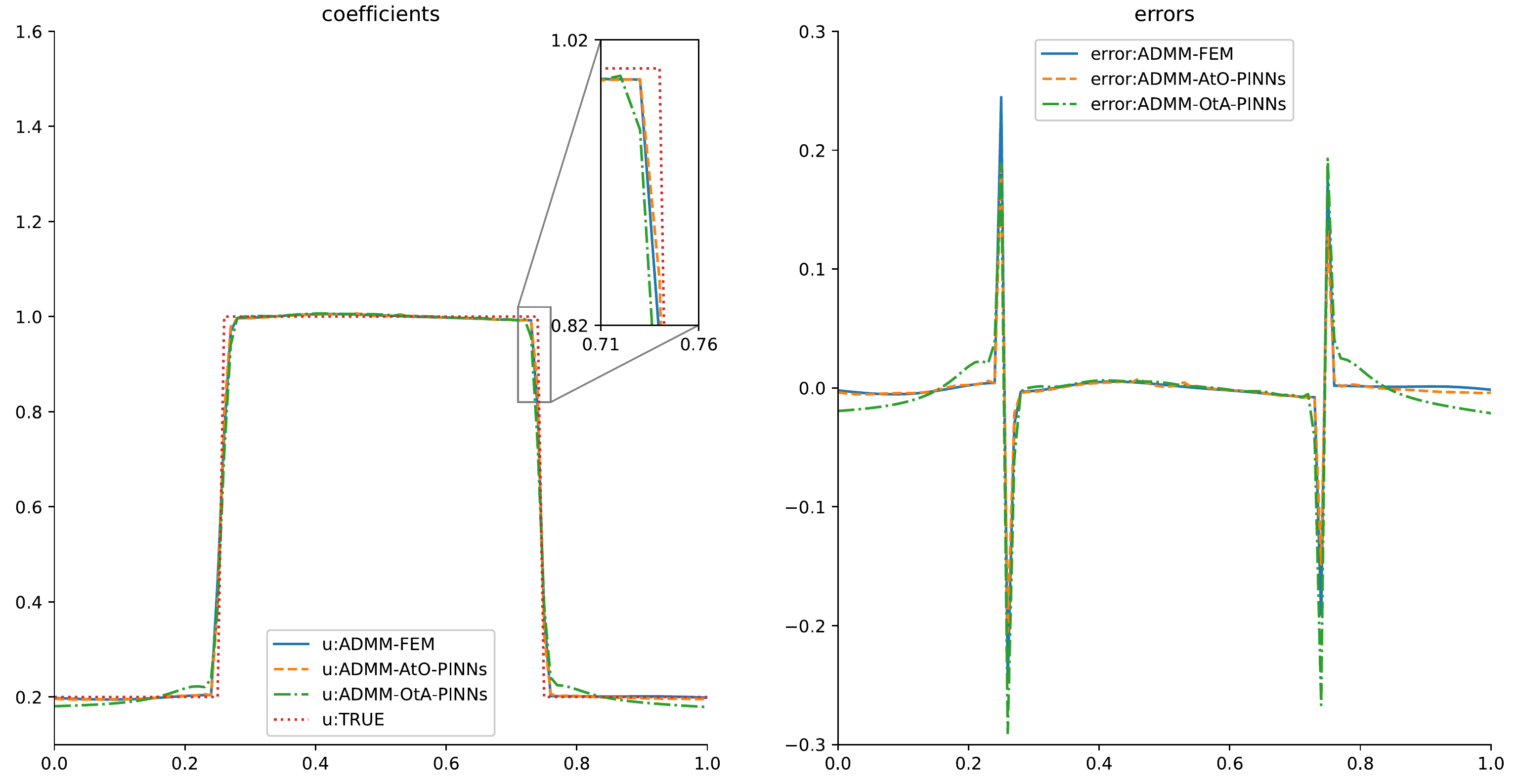}}
	
	\subfloat[Recovered states  and errors by different algorithms]{\includegraphics[width=\textwidth]{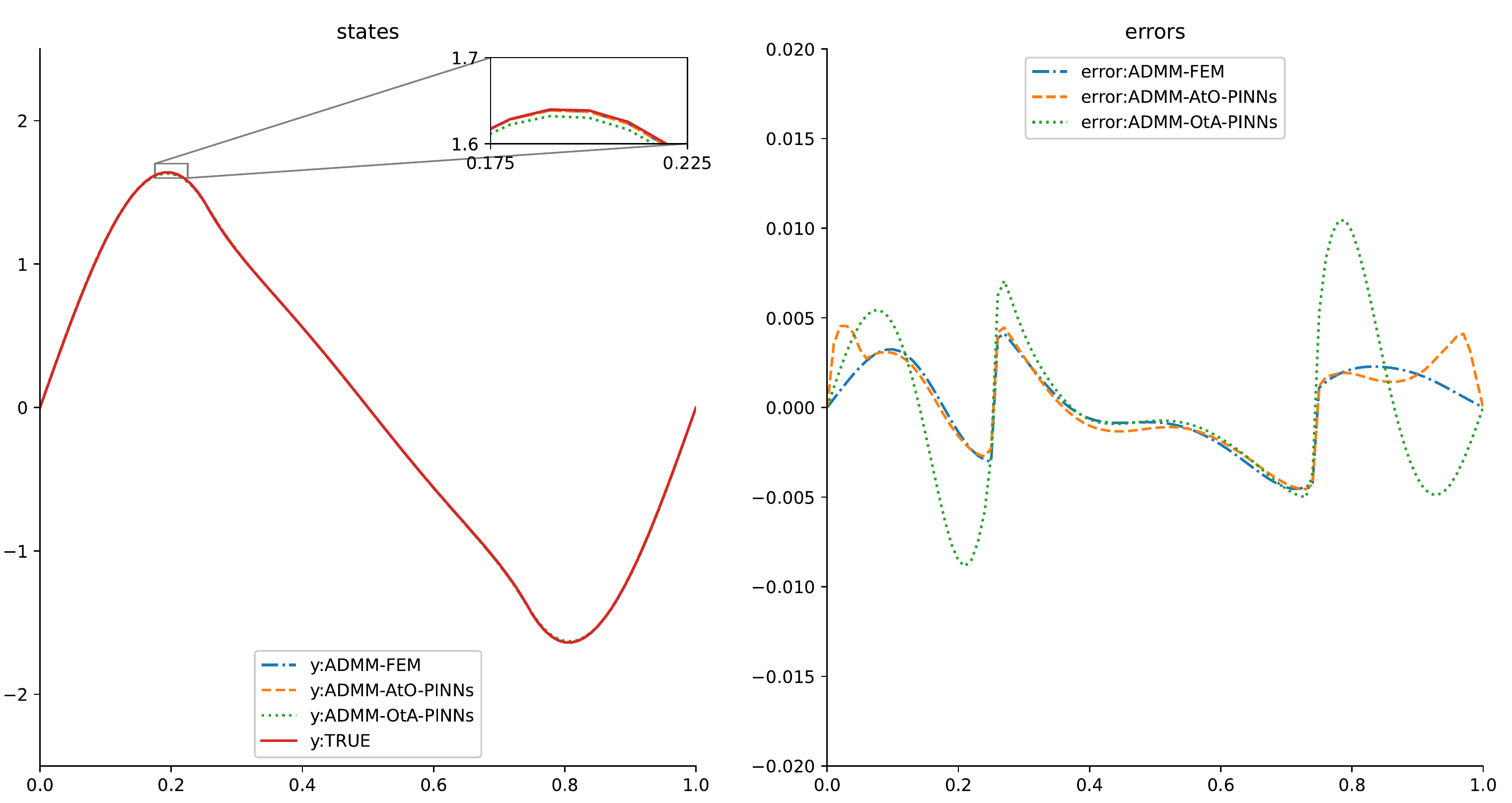}}
	
	\centering
\end{figure}

In the implementation of the  ADMM-AtO-PINNs, the ADMM-OtA-PINNs and the ADMM-Gauss-Newton, we execute 20 outer ADMM iterations with $\beta=0.1, z^0=0, \lambda^0=0$. The numerical results are reported in Table \ref{tab:err_ex1} and Figure \ref{fig:reaction_id}. From these results, we can see that the $y$ and $u$ obtained by the ADMM-AtO-PINNs and the ADMM-OtA-PINNs are as accurate as those by the ADMM-Gauss-Newton, and all of them are good approximations of the exact ones.
\begin{table}[h]
	\setlength{\abovecaptionskip}{0pt}
	\setlength{\belowcaptionskip}{3pt}
	\centering
	\caption{Numerical errors of different algorithms for Example 1. }\label{tab:err_ex1}
	{\small\begin{tabular}{|c|c|c|c|c|}
			\hline
			Algorithm&ADMM-AtO-PINNs&ADMM-OtA-PINNs&ADMM-Gauss-Newton\\
			\hline
			$\frac{\|u-u_{true}\|}{\|u_{true}\|}$&0.0472 &0.0674 &0.0609\\
			\hline
		\end{tabular}
	}
\end{table}

\section{Control Constrained Optimal Control of the Burgers Equation}\label{se:control_burgers}

As a simplified model for convection-diffusion phenomena, the Burgers
equation plays a crucial role in diverse physical problems such as shock
waves, supersonic flow, traffic flows, and acoustic transmission. Its study usually serves
as a first approximation to more complex convection diffusion phenomena. Hence, numerical study for
the optimal control of the Burgers equation is important for the development of numerical methods for optimal control of more complicated models in fluid dynamics like the Navier-Stokes equations. We refer to e.g., \cite{volkwein1997,volkwein2001} for more discussions.

\subsection{Optimal Control Model} Let $\Omega=(0,1)$, we consider the control constrained optimal control of the
stationary Burgers equation as that in \cite{DeK2004}:
 \begin{equation}\label{model_burgers}
\left\{
\begin{aligned}
	\min_{y\in Y, u\in U} &\frac{1}{2}\|y - y_d\|^2_{L^2(\Omega)} + \frac{\alpha}{2}\|u\|^2_{L^2(\omega)}+I_{U_{ad}}(u)\\
	\text{s.t.}& -\nu y''+yy'=\chi_{\omega}u \quad \text{in}~(0,1),\qquad y(0)=y(1)=0.
\end{aligned}
\right.
 \end{equation}
Above, $Y=H_0^1(\Omega)$, $U=L^2(\omega)$, the regularization parameter $\alpha>0 $, and the parameter $\nu > 0$ denotes the viscosity coefficient of the fluid, which is equal to $\frac{1}{Re}$ with $Re$ the Reynolds number. The target function $y_d\in H^1(\Omega)$ is given, $\chi_{\omega}$ is the characteristic function of $\omega\subset \Omega$ defined by $\chi_{\omega}(x)=1$ if $x\in \omega$ and $\chi_{\omega}(x)=0$ otherwise, and $I_{U_{ad}}$ is the indicator function of the control set $U_{ad}$ defined by
 $$
U_{ad}=\{u\in L^2(\Omega)| a\leq u(x)\leq b,~\text{a.e.~in}~\omega\},
 $$
where $a$ and $b$ are given constants. The optimal control of the Burgers equation has been extensively studied in the literature, see e.g.,  \cite{DeK2004,volkwein1997,volkwein2001} and references therein. In particular, a complete analysis
of the optimal control of the Burgers equation was given in \cite{volkwein1997} and some numerical studies can be referred to \cite{DeK2004}.

\subsection{ADMM-PINNs for (\ref{model_burgers})}
Let $y(u)$ be the solution of the state equation (\ref{state_equation}) corresponding to $u$. By introducing $z\in L^2(\omega)$ satisfying $u=z$, problem (\ref{model_burgers}) can be reformulated as
 \begin{equation}\label{model_burgers_e}
	\begin{aligned}
		\underset{u,z\in L^2(\omega) }{\min}~& \frac{1}{2}\|y(u)-y_d\|^2_{L^2(\Omega)}+\frac{\alpha}{2}\|u\|^2_{L^2(\omega)}+I_{U_{ad}}(z),\quad
		\text{s.t.}~&~u=z.
	\end{aligned}
 \end{equation}
The augmented Lagrangian functional associated with (\ref{model_burgers_e}) is defined as
 \begin{equation}\label{L:burgers_e}
L_{\beta}(u,z;\lambda)=\frac{1}{2}\|y(u)-y_d\|^2_{L^2(\Omega)}+\frac{\alpha}{2}\|u\|^2_{L^2(\omega)}+I_{U_{ad}}(z)-(\lambda,u-z)_{L^2(\omega)}+\frac{\beta}{2}\|u-z\|_{L^2(\omega)}^2.
 \end{equation}
 The ADMM algorithm for (\ref{model_burgers}) follows from \eqref{admm} with $L_{\beta}(u,z;\lambda)$ defined by \eqref{L:burgers_e}.

It is easy to show that the $z$-subproblem (\ref{admm_z}) reads as
 $$
z^{k+1}=\arg\min_{z\in L^2(\omega)}I_{U_{ad}}(z)-(\lambda^k,u^{k+1}-z)_{L^2(\omega)}+\frac{\beta}{2}\|u^{k+1}-z\|_{L^2(\omega)}^2,
 $$
which  has a closed-form solution
$
z^{k+1}=\mathbb{P}_{U_{ad}}(u^{k+1}-\frac{\lambda^k}{\beta}).
$
Moreover, the $u$-subproblem (\ref{admm_u}) can be reformulated as the following unconstrained optimal control problem
{\small \begin{equation*}
	\min_{y\in L^2(\Omega),u\in L^2(\omega)}\mathcal{J}^k_{BC}(y,u):=\frac{1}{2}\|y-y_d\|^2_{L^2(\Omega)}+\frac{\alpha}{2}\|u\|^2_{L^2(\omega)}-(\lambda^k,u-z^k)_{L^2(\omega)}+\frac{\beta}{2}\|u-z^k\|_{L^2(\omega)}^2,
\end{equation*}}
subject to the Burgers equation in (\ref{model_burgers}). Also, if $u$ is a local optimal solution of problem (\ref{admm_u}), and $y$ and $p$ are the corresponding state and adjoint variables, then they satisfy the following first-order optimality system:
 \begin{equation}\label{oc_u_burgers}
	\left\{
	\begin{aligned}
		&p|_{\omega} +(\alpha+\beta) u-\lambda^k-\beta z^k=0,\\
		& -\nu y''+yy'=\chi_{\omega}u~ \text{in}~(0,1),~ y(0)=y(1)=0,\\
		&-\nu p''-yp'= y-y_d~ \text{in}~(0,1),~p(0)=p(1)=0.
	\end{aligned}
	\right.
 \end{equation}
Using the relation $u=\frac{1}{\alpha+\beta}(-p|_{\omega}+\lambda^k+\beta z^k)$, we can remove the control $u$ in (\ref{oc_u_burgers}) and obtain the following reduced optimality system:
 \begin{equation}\label{oc_u_burgers_red}
	\left\{
	\begin{aligned}
		& -\nu y''+yy'=\chi_{\omega}(\frac{1}{\alpha+\beta}(-p|_{\omega}+\lambda^k+\beta z^k))~ \text{in}~(0,1),~ y(0)=y(1)=0,\\
		&-\nu p''-yp'= y-y_d ~ \text{in}~(0,1),~p(0)=p(1)=0.
	\end{aligned}
	\right.
 \end{equation}

As discussed in Section \ref{se:solution_u},  an ADMM-AtO-PINNs and an ADMM-OtA-PINNs can be specified from Algorithm \ref{alg:admm_pinn} for (\ref{model_burgers}). Note that he ADMM-OtA-PINNs is implemented with the reduced optimality system (\ref{oc_u_burgers_red}), instead of the full optimality system (\ref{oc_u_burgers}).
\subsection{Numerical Results for  (\ref{model_burgers})}

In this subsection, we test the control constrained optimal control of the Burgers equation (\ref{model_burgers}) and numerically verify the effectiveness of the proposed  ADMM-AtO-PINNs and ADMM-OtA-PINNs.

\medskip

\noindent\textbf{Example 2.} We consider the example given in \cite{DeK2004}. In particular, we set $\Omega=(0,1)$, $\omega=\Omega$, and $\nu = \frac{1}{12}, \alpha = 0.1,  a=-\infty, b =
0.3$, and $y_d= 0.3$ in problem (\ref{model_burgers}).

For the implementation of PINNs,  we approximate $u$ by a fully-connected feed-forward neural network $\hat{u}(x,\bm{\theta}_u)$. Moreover, $y$ and $p$ are respectively approximated by $\hat{y}(x;\bm{\theta}_y)=x(x-1)\mathcal{NN}(x;\bm{\theta}_y)$ and $\hat{p}(x;\bm{\theta}_p=x(x-1)\mathcal{NN}(x;\bm{\theta}_p)$, with $\mathcal{NN}(x;\bm{\theta}_y)$ and $\mathcal{NN}(x;\bm{\theta}_p)$ being fully-connected feed-forward neural networks, so that the boundary conditions  $y(0)=y(1)=0$ and $p(0)=p(1)=0$ can be satisfied automatically. Hence, we set $w_b=0$. All other weights in the loss function are set to be 1.  All the neural networks consist of 4 hidden layers with 50 neurons per hidden layer and hyperbolic tangent activation functions.  To evaluate the loss, we take the training sets $\mathcal{T}_i\subset\Omega$ as the grid nodes of the lattice generated by uniformly partitioning $\Omega=(0,1)$ with mesh size $h=1/100$.  To train the neural networks,  we first implement 5000 iterations of the Adam \cite{kingma2015} with learning rate $\eta=10^{-3}$ to obtain good guesses of the parameters and then switch to the L-BFGS \cite{nocedal1980} with a strong Wolfe line search strategy for 1000 iterations to achieve higher accuracy.  The neural network parameters are initialized using the
default initializer of Pytorch.   In the implementation of the ADMM-AtO-PINNs and ADMM-OtA-PINNs, we execute 20 outer ADMM iterations with $\beta=0.1, z^0=0$, and $\lambda^0=0$.

Moreover, to validate the effectiveness of the ADMM-AtO-PINNs and the ADMM-OtA-PINNs, we compare our results with the reference ones obtained by high-fidelity traditional numerical methods based on FEM. To be concrete, problem (\ref{model_burgers}) is first discretized by FEM and then is solved by the SSN method  \cite{DeK2004}. A FEM-SSN method is thus specified for solving (\ref{model_burgers}).  To implement the FEM, the domain $\Omega=(0,1)$ is uniformly partitioned with mesh size $h=1/100$.

\begin{figure}[htpb]
	\caption{\small Numerical results for Example 2}\label{fig:control_burgers}
	\centering
{\includegraphics[width=\textwidth]{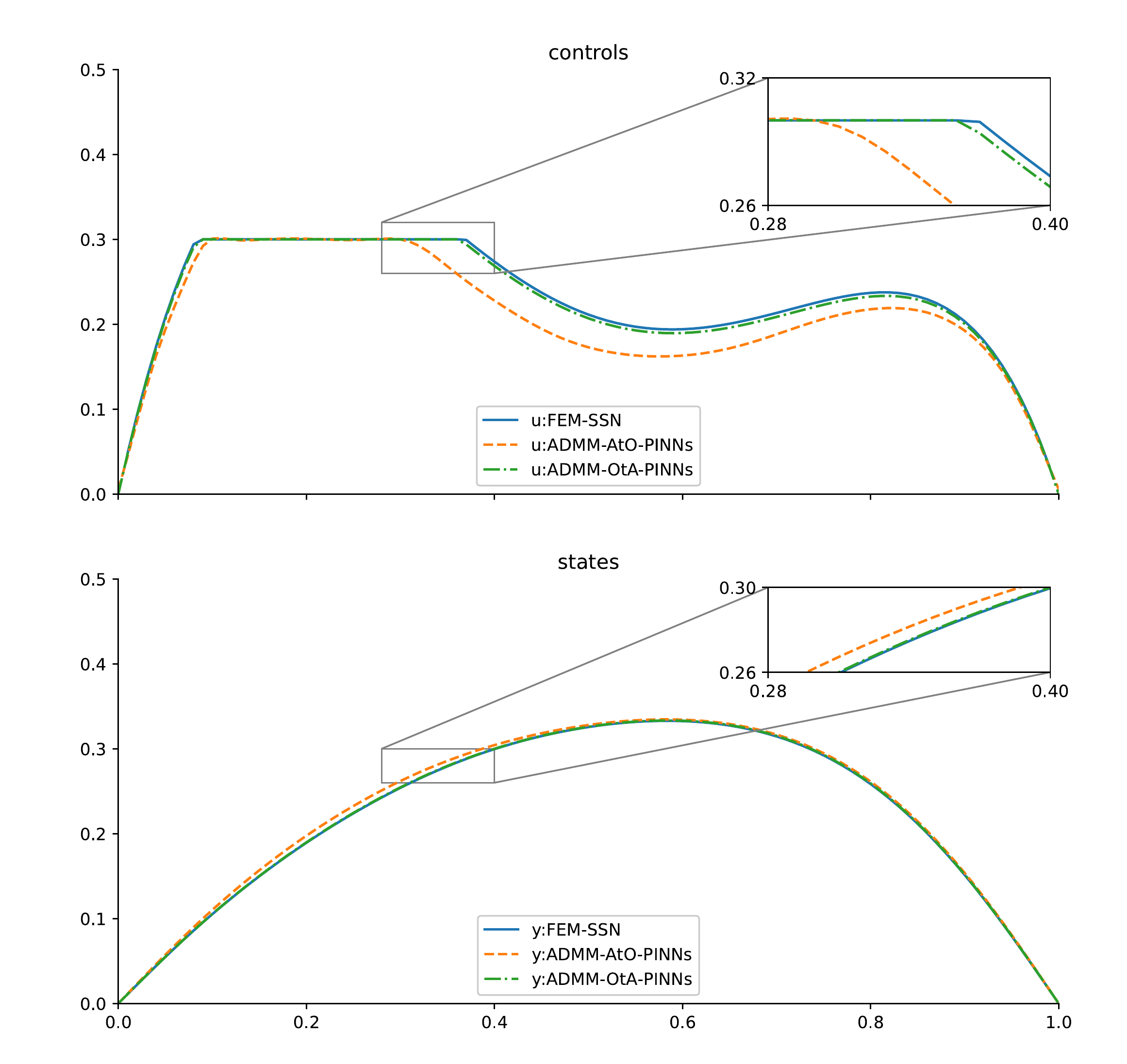}}
\end{figure}

Figure \ref{fig:control_burgers} shows the computed optimal controls and states by different algorithms.  We observe that the states obtained by all the algorithms are almost the same. The control constraint $u(x)\leq 0.3$ a.e. in $\Omega$ is satisfied in all cases. Furthermore, the control computed by  the ADMM-OtA-PINNs is so close to the one by the FEM-SSN and they cannot be visually distinguished. However, we note that the control computed by  the ADMM-AtO-PINNs is different, see Table \ref{tab:err_ex2} for further verification. A possible reason is that problem (\ref{model_burgers}) has a non-unique solution due to its nonconvexity, and  the ADMM-AtO-PINNs finds a different control that solves the problem. All these results validate the effectiveness of our proposed ADMM-PINNs algorithms.

\begin{table}[h]
	\setlength{\abovecaptionskip}{0pt}
	\setlength{\belowcaptionskip}{3pt}
	\centering
	\caption{Relative difference between FEM-SSN and ADMM-PINNs for Example 2.}\label{tab:err_ex2}
	{\small\begin{tabular}{|c|c|c|c|}
			\hline
			Algorithm&ADMM-AtO-PINNs&ADMM-OtA-PINNs\\
			\hline
			$\frac{\|u-u_{FEM}\|}{\|u_{FEM}\|}$&0.101 &0.018 \\
			\hline
		\end{tabular}
	}
\end{table}

\section{Discontinuous Source Identification for Elliptic PDEs}\label{se:source_iden}
In this section,  we discuss the implementation of the ADMM-PINNs algorithmic framework to discontinuous source identification for elliptic PDEs.
\subsection{Problem Statement}
We consider the problem of identifying a discontinuous source in a prototypical elliptic PDE from Dirichlet boundary observation data:
 \begin{equation}\label{model:si}
	\begin{aligned}
		\min_{u\in U,y\in Y}&~\Big\{\frac{1}{2}\int_{\partial\Omega} | y- y^{\delta}|^2 dx+\gamma\int_{\omega}|\nabla u|dx\Big\},\\
		\text{s.t.}~&-\Delta y+cy=u~\text{in}~\Omega,~\frac{\partial y}{\partial n}=0~\text{on}~\partial\Omega,~ u=0~\text{in}~\Omega/\omega
	\end{aligned}
 \end{equation}
where $\Omega$ is a bounded domain in $\mathbb{R}^d$ ($d\geq 1$) with a piecewise smooth boundary $\partial\Omega$, $\omega\subset \Omega$ is an open set, $Y=L^2(\Omega)$, and $U=L^2(\Omega)\cap BV(\Omega)$ is a Banach space endowed with the norm $\|u\|_{U}=\|u\|_{L^2(\Omega)}+\|u\|_{BV(\Omega)}$ \cite{chavent1997regularization}. The constant coefficient $c>0$ is given, $u$ is the unknown source, and $n$ denotes the outwards pointing unit normal
vector of the boundary $\partial\Omega$.   Here, we consider the case that some observation data, denoted by $y^{\delta}\in L^2(\Omega)$, of the solution $y$ on the boundary $\partial\Omega$ is available subject to some noise with the noisy level $\delta>0$.  In other words, we attempt to use the Dirichlet boundary observation data $\ y^{\delta}$ to identify the unknown source $u$ by solving (\ref{model:si}). Moreover,  $\frac{1}{2}\int_{\partial\Omega} | y- y^{\delta}|^2\mathrm{d}x$ is a data-fidelity term, $\int_{\omega}|\nabla u|dx$ is the TV-regularization defined by (\ref{eq:TV}), and $\gamma> 0$ is a regularization parameter. Problem (\ref{model:si}) and its variants have found various applications in crack determination \cite{alves2004}, electroencephalography \cite{baillet2001} and inverse electrocardiography \cite{wang2013}.

\subsection{ADMM-PINNs for (\ref{model:si})}

For any given $u\in L^2(\Omega)$, let $y(u)$ be the solution of the equation
 $$
-\Delta y+cy=u~\text{in}~\Omega,~\frac{\partial y}{\partial n}=0~\text{on}~\partial\Omega.
 $$
 By introducing an auxiliary variable $z\in U$ such that $u=z$, we can reformulate problem (\ref{model:si}) as
  \begin{equation}\label{model:sie}
 	\begin{aligned}
 		\min_{u,z\in U}~\frac{1}{2}\int_{\partial\Omega} | y(u)- y^{\delta}|^2\mathrm{d}x+\gamma\int_{\omega}|\nabla z|dx,\quad
 		\text{s.t.}~u=z ~\text{in}~\Omega,~z=0 ~\text{in}~\Omega/\omega
 	\end{aligned}
  \end{equation}
Since $U= L^2(\Omega)\cap BV(\Omega)$, the $L^2$-inner product and the $L^2$-norm are well-defined in $U$. Then, the augmented Lagrangian functional associated with (\ref{model:sie}) is defined as
  \begin{equation}\label{L:sie}
L_{\beta}(u,z;\lambda)= \frac{1}{2}\int_{\partial\Omega} | y(u)- y^{\delta}|^2\mathrm{d}x+\gamma\int_{\omega}|\nabla z|dx-(\lambda,u-z)_{L^2(\Omega)}+\frac{\beta}{2}\|u-z\|_{L^2(\Omega)}^2,
  \end{equation}
where $\beta>0$ is a penalty parameter. Then the ADMM algorithm for (\ref{model:si}) follows from \eqref{admm} with $L_{\beta}(u,z;\lambda)$ defined by \eqref{L:sie}.

We first note that the $z$-subproblem (\ref{admm_z}) is a TV-regularized optimization problem
 \begin{equation}\label{admm_z_e}
z^{k+1}=\arg\min_{z\in U} \gamma\int_{\omega}|\nabla z|dx+\frac{\beta}{2}\|z-(u^{k+1}-\frac{\lambda^k}{\beta})\|^2,\quad\text{in}~\omega,\\
 \end{equation}
and $z^{k+1}=0$ in ${\Omega/\omega}$.
We further assume that $\omega\subset \mathbb{R}^2$ is a rectangular domain. Then, inspired by \cite{tyy2021}, problem \eqref{admm_z_e} after some proper discretization can be solved by a pre-trained deep CNN. To elaborate,
let $R(z):=\int_{\omega}|\nabla z|dx$, then the proximal operator of $\gamma R(z)$ is given by
\begin{equation*}
  \text{Prox}_{\gamma R(z)}(v)=\arg \min_{z\in U} \Big\{\gamma R(z)+\frac{1}{2}\|z-v\|_{L^2(\omega)}^2\Big\},~\forall v\in L^2(\omega).
\end{equation*}
Then, the solution of \eqref{admm_z_e} can be expressed as
 $$
z^{k+1}=\text{Prox}_{\frac{\gamma}{\beta}R(z)}(u^{k+1}-\frac{\lambda^k}{\beta})\quad\text{in}~\omega.
 $$
Following  \cite{RudinOF:1992}, $\text{Prox}_{\frac{\gamma}{\beta} R(z)}$ can be interpreted as an image denoising operator.

Suppose that the domain $\omega$ is uniformly partitioned with mesh size $h$, and $\bm{x}=\{x_i\}_{i=1}^N$ are the grid nodes of the resulting lattice. For any function $q:\omega\subset\mathbb{R}^2\rightarrow \mathbb{R}$, we let  $q_h:\mathbb{R}^2\rightarrow \mathbb{R}$ defined on the lattice be a finite approximation of $q$. Then, there exists a one-to-one mapping between $q_h$ and an $m\times n$ raster image $\mathcal{I}$ ($m\times n=N$), where the gray value at pixel $x_i$ of the image $\mathcal{I}$ corresponds to the value of the function $q_h$ at node $x_i$.  Then, pre-trained deep CNNs, which have been widely used for various image denoising problems, can be applied. For this purpose,  let $\mathcal{M}$ denote the mapping from $q_h$ to the raster image $\mathcal{I}$, and $\mathcal{C}_{\sigma}$ represents a pre-trained deep CNN with $\sigma$ the variance of the noise used for training the CNN. Solving problem \eqref{admm_z_e}  by a pre-trained deep CNN is summarized in Algorithm \ref{Algo_7}. Note that once a pre-trained deep CNN is applied, the parameter $\gamma$ is not needed any more and its role will be replaced by the denoising parameter $\sigma$.
 \begin{algorithm}[htpb]\caption{A deep CNN-based method for  (\ref{admm_z_e}).}\label{Algo_7}
  \begin{algorithmic}[1]
      \STATE $\mathcal{I}^{k+1}_{input}:=\mathcal{M}(u_h^{k+1}(\bm{x})-\frac{\lambda_h^k(\bm{x})}{\beta})$.
      \STATE $\mathcal{I}^{k+1}_{output}:=\mathcal{C}_{\sigma }(\mathcal{I}_{input}^{k+1})$.
      \STATE $z_h^{k+1}(\bm{x})=\mathcal{M}^{-1}(\mathcal{I}^{k+1}_{output})$.
  \end{algorithmic}
\end{algorithm}

The $u$-subproblem (\ref{admm_u}) can be reformulated as the following PDE-constrained optimization problem:
 $$
\min \mathcal{J}^k_{SI}(y,u):=\frac{1}{2}\int_{\partial\Omega} | y- y^{\delta}|^2\mathrm{d}x-(\lambda^k, u-z^k)+\frac{\beta}{2}\|u-z^k\|^2,
 $$
subject to the equation in (\ref{model:si}).  Its first-order optimality system reads as:
 \begin{equation}\label{oc_u_SI}
	\left\{
	\begin{aligned}
		&p+\beta u-\lambda^k-\beta z^k=0,\\
		&-\Delta y+cy=u \text{~in~} \Omega,~ \frac{\partial y}{\partial n}=0 \text{~on~} \partial \Omega,\\
		&-\Delta p+cp=0 \text{~in~} \Omega,~ \frac{\partial p}{\partial n}=y-y^{\delta} \text{~on~} \partial \Omega,
	\end{aligned}
	\right.
 \end{equation}
where $u$ is the solution to the $u$-subproblem (\ref{admm_u}), $y$ is the corresponding solution of the equation in (\ref{model:si}), and $p$ is the adjoint variable.
Furthermore, using the relation $u=\frac{-p+\lambda^k+\beta z^k}{\beta}$, we can rewrite \eqref{oc_u_SI} as
 \begin{equation}\label{oc_u_SI_r}
	\left\{
	\begin{aligned}
		&-\Delta y+cy=\frac{-p+\lambda^k+\beta z^k}{\beta} \text{~in~} \Omega,~ \frac{\partial y}{\partial n}=0 \text{~on~} \partial \Omega,\\
		&-\Delta p+cp=0 \text{~in~} \Omega,~ \frac{\partial p}{\partial n}=y-y^{\delta} \text{~on~} \partial \Omega.
	\end{aligned}
	\right.
 \end{equation}

As discussed in Section \ref{se:solution_u}, the $u$-subproblem (\ref{admm_u}) can  be solved by applying Algorithm \ref{alg:pinn_d} or  Algorithm \ref{alg:pinn_oc}. Accordingly, we can obtain two ADMM-PINNs algorithms for solving problem (\ref{model:si}). However,  it was empirically observed that, when the ADMM-AtO-PINNs is applied, the corresponding neural networks are hard to train and the numerical results are not satisfactory. A main reason is that the balance between the objective function and the PDE losses is sensitive and unstable for this case. To address this issue, much more effort is needed to find suitable loss weights, which is beyond the scope of the paper. Hence, we only focus on the application of the ADMM-OtA-PINNs to problem (\ref{model:si}).  Note that the ADMM-OtA-PINNs is implemented with the reduced optimality system (\ref{oc_u_SI_r}), instead of the full optimality system (\ref{oc_u_SI}).

\subsection{Numerical Results for (\ref{model:si})}

In this section, we devote ourselves to showing the numerical results of the ADMM-OtA-PINNs for solving  (\ref{model:si}).

\noindent\textbf{Example 3.}
We set $\Omega=(0,1)\times(0,1)$, $\omega=\omega_1\cup\omega_2$ with $\omega_1=(0.25, 0.5)\times(0.25,0.75), \omega_2=(0.5,0.75)\times(0.25,0.75)$ and $c=1$.
We set the true source $u$ as
 $$
u_{true}(x)=3\chi_{w_1}(x)-9\chi_{\omega_2}(x),
 $$
where $\chi_{w_1}$ and $\chi_{w_2}$ are the characteristic functions of $\omega_1$ and $\omega_2$, respectively.
The observed data $y^\delta=y_{true}+\delta \| y_{true}\| \text{rand} (x)$, where $\delta=0.1$, $y_{true}$ is the solution of equation \eqref{model:si} corresponding to $u_{true}$, and $\text{rand} (x)$ is a standard normal distributed random function.   In our numerical experiments, $y_{true}$ is the solution of \eqref{eq:irc} obtained from $u_{true}$ by FEM. To implement the FEM used for this example, the domain $\Omega$ is uniformly triangulated with mesh size $h=1/64$.

For the implementation of PINNs,  we approximate $y$ and $p$ by $\hat{y}(x;\bm{\theta}_y)=\mathcal{NN}(x;\bm{\theta}_y)$ and $\hat{p}(x;\bm{\theta}_p)=\beta\mathcal{NN}(x;\bm{\theta}_p)$, with $\mathcal{NN}(x;\bm{\theta}_y)$ and $\mathcal{NN}(x;\bm{\theta}_p)$ being fully-connected feed-forward neural networks.  All the neural networks consist of 3 hidden layers with 32 neurons per hidden layer and hyperbolic tangent activation functions.  To evaluate the loss, we take the training sets $\mathcal{T}_i\subset\Omega$ as $31\times 31$ uniformly sampled points and $\mathcal{T}_b\subset\partial\Omega$ as $4\times 128$ uniformly sampled points. The weights $\frac{w_i}{|\mathcal{T}_i|}$,  $\frac{w_b}{|\mathcal{T}_b|}$,$w_o$, $w_e$, $w_y$, and $w_p$ are set to be 1. To train the neural networks,  we first execute the Adam \cite{kingma2015} with learning rate $\eta=10^{-4}$ for  20000 iterations to obtain good guesses of the neural network parameters and then switch to the L-BFGS \cite{nocedal1980} with a strong Wolfe line search strategy for 1000 iterations to achieve higher accuracy.  The neural network parameters are initialized using the default initializer of Pytorch.

For the CNN used in this example, we use the same network architecture, the same training dataset, and the PyTorch version of the original code as \cite{ZhangZCMZ:2017} (see \url{https://github.com/cszn/DnCNN} for the original code and \url{https://github.com/cszn/KAIR} for the PyTorch version) to train 30 CNNs for the cases where $\sigma=1,2,\cdots, 30$. Hence, $C_\sigma$ in Algorithm \ref{Algo_7} is chosen as one of these pre-trained 30 CNNs with a specified value of $\sigma$. According to our tests, $\sigma=24$ is the best choice. The mapping $\mathcal{M}$ in Algorithm \ref{Algo_7} is specified as $\mathcal{M}(z_h)=(\bm{z}+10)/20\times256$, where $\bm{z}=z_h(\bm{x})$ and $\mathcal{I}$ is a raster image. To avoid the gray values of the raster image exceeding the range, we added a projection to ensure the gray values of $\mathcal{I}$ strictly fall in $(0,255)$.

Moreover, to validate the effectiveness of the ADMM-OtA-PINNs, we compare our results with the reference ones obtained by a high-fidelity FEM-based numerical method. To be concrete, the subproblem (\ref{admm_u}) is first discretized by FEM and then is solved by a conjugate gradient (CG) method.  An ADMM-CG-FEM algorithm is thus specified for solving (\ref{model:si}).  Similar methods can be found in \cite{GSY2019,GSYY2022}. For implementing the ADMM-OtA-PINNs and the ADMM-CG-FEM, we execute 350 ADMM iterations with $\beta=5\times10^{-2}, z^0=0, \lambda^0=0$.

\begin{figure}[htpb]
	\caption{\small Numerical results for Example 3. Noise level $\delta=10\%.$ Relative error of $u$: 0.174639 (ADMM-CG), 0.171283 (ADMM-OtA-PINNs) }\label{fig:SI}
	\centering
	\subfloat[True source]{\includegraphics[width=0.32\textwidth]{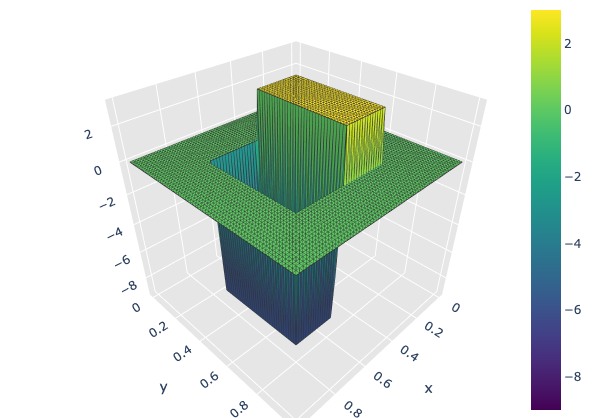}}
	\subfloat[Recovered source by ADMM-CG-FEM]{\includegraphics[width=0.32\textwidth]{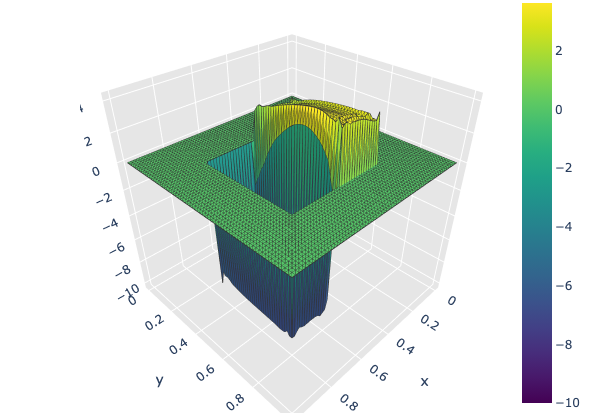}}
		\subfloat[Recovered source by  ADMM-OtA-PINNs]{\includegraphics[width=0.32\textwidth]{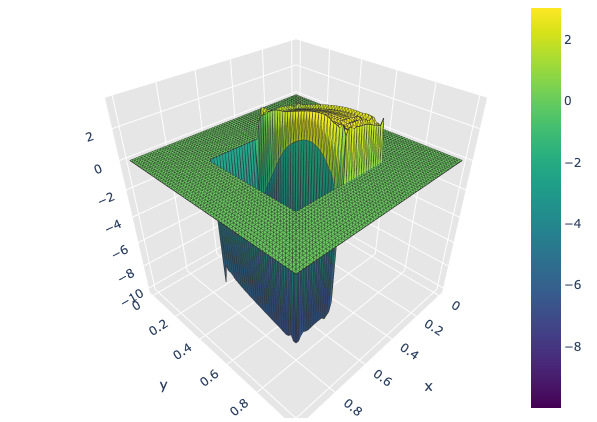}}
		
				\subfloat[True source]{\includegraphics[width=0.32\textwidth]{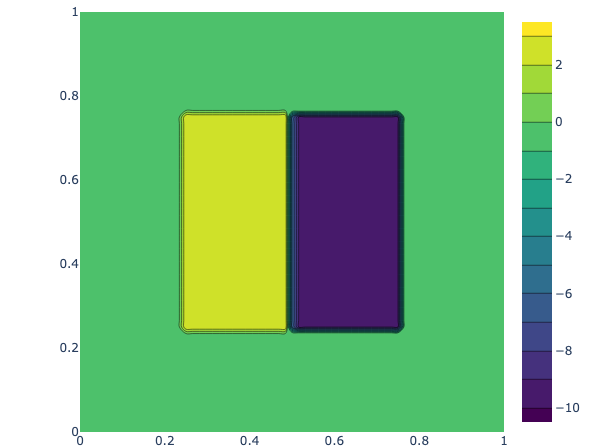}}				
			\subfloat[Recovered source by ADMM-CG-FEM]{\includegraphics[width=0.32\textwidth]{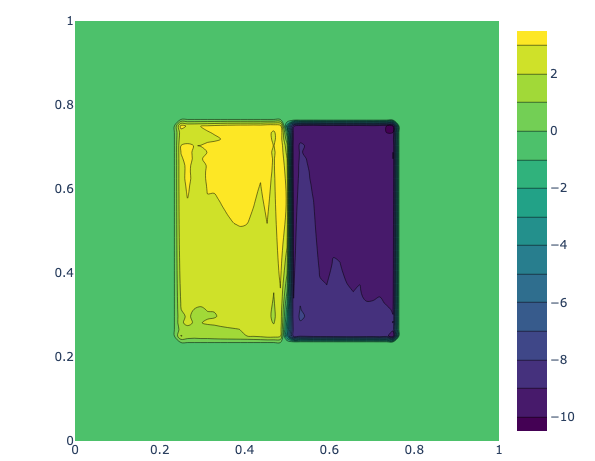}}
			\subfloat[Recovered source by  ADMM-OtA-PINNs]{\includegraphics[width=0.32\textwidth]{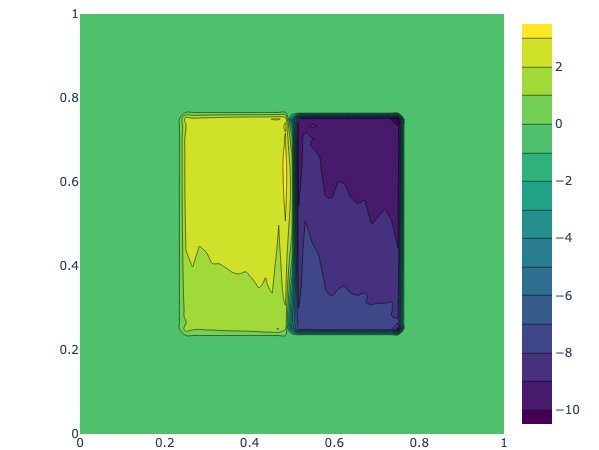}}	
			
			\subfloat[True state]{\includegraphics[width=0.32\textwidth]{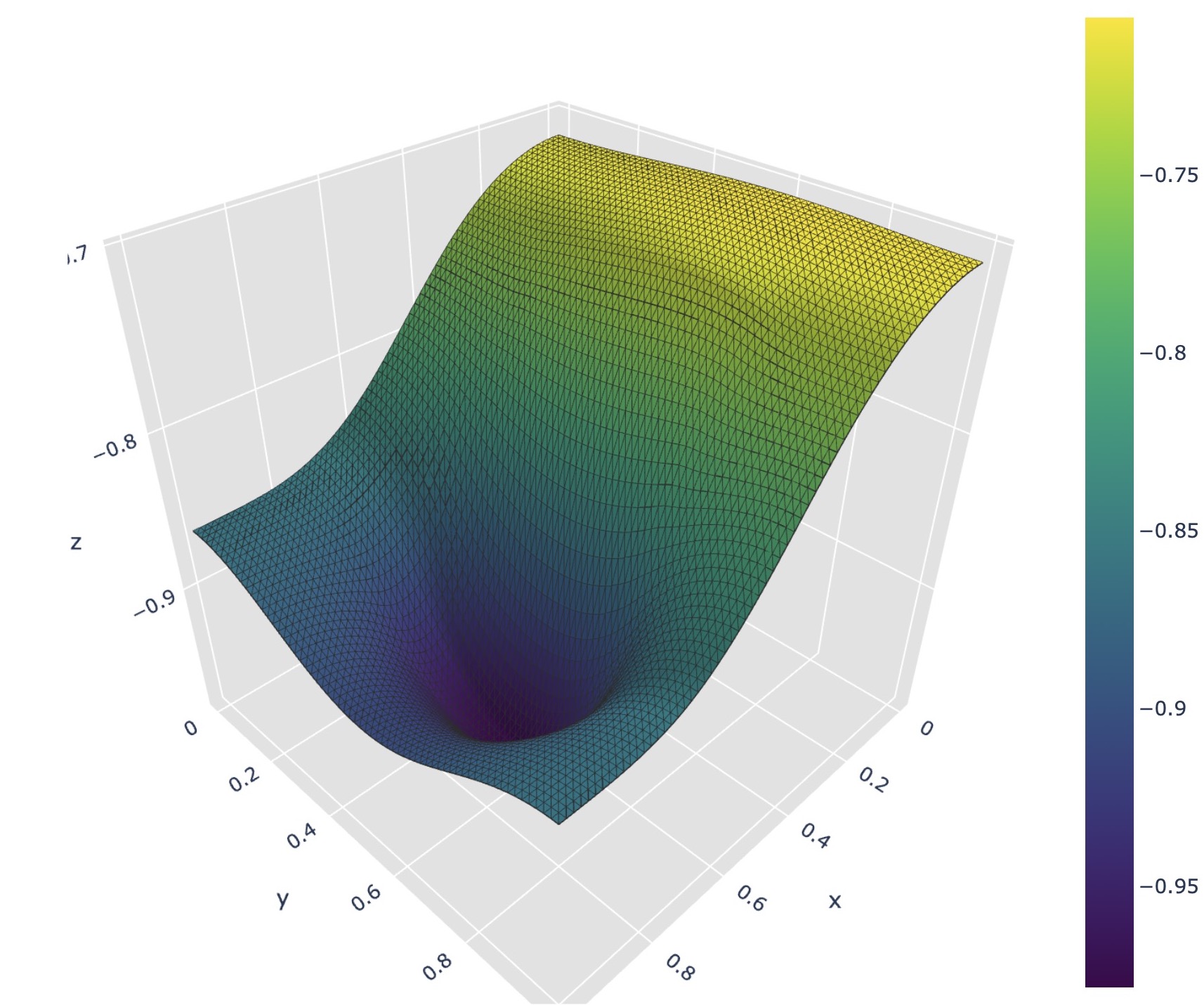}}		
	\subfloat[Recovered state  by ADMM-CG-FEM]{\includegraphics[width=0.32\textwidth]{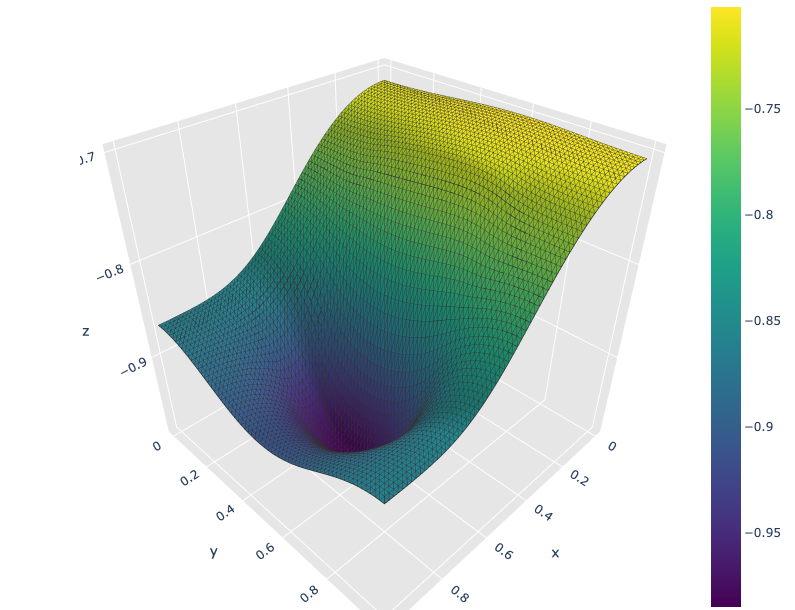}}
		\subfloat[Recovered state  by  ADMM-OtA-PINNs]{\includegraphics[width=0.32\textwidth]{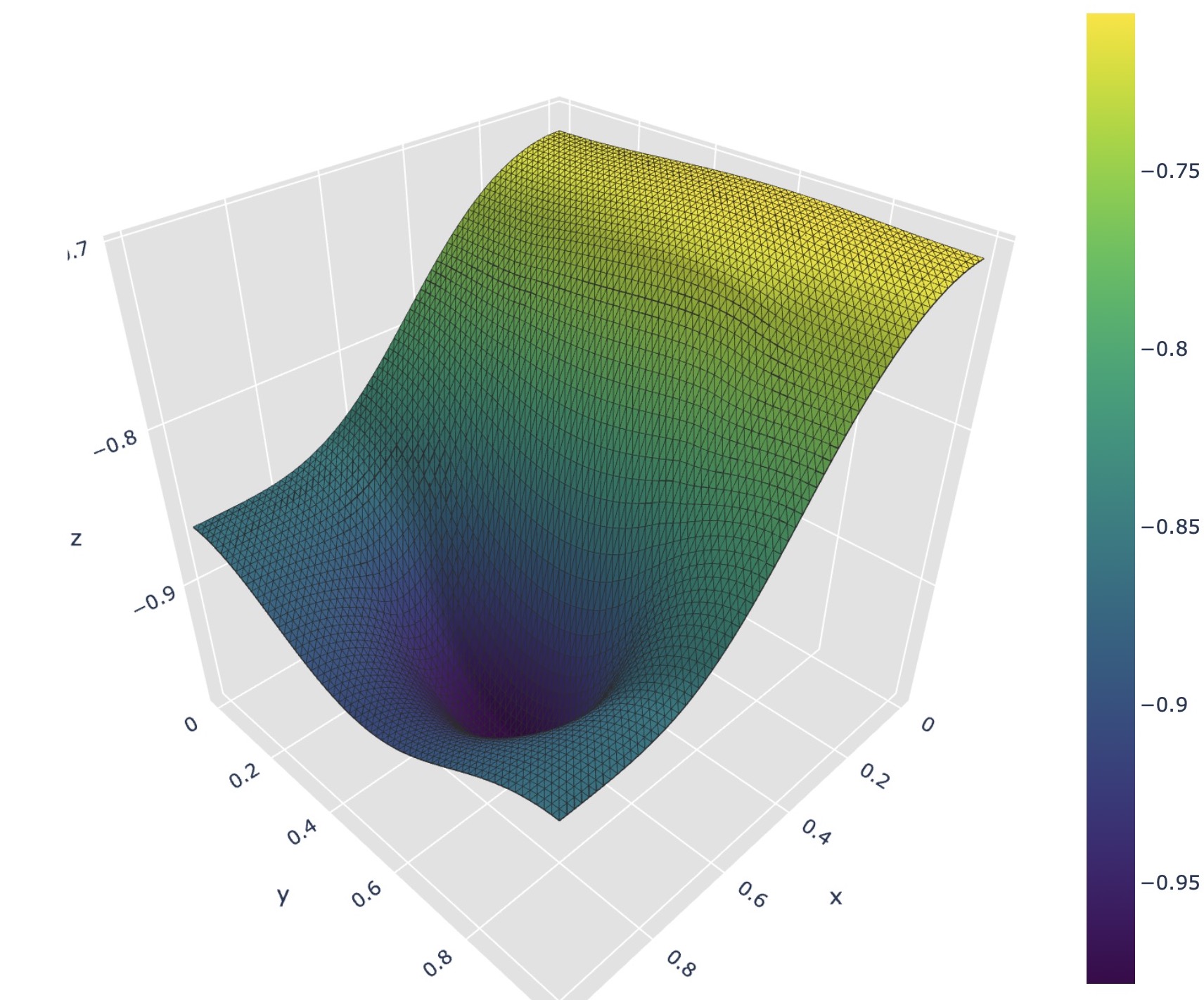}}
	\centering
\end{figure}

We observe from the results that the ADMM-OtA-PINNs is robust to the noise level $\delta$. We report the recovered sources and states by different algorithms with $\delta=10\%$ in
Figure \ref{fig:SI}. We can see that even if the noise increases to $10\%$, the identified
solutions are still very satisfactory. In particular, the computed states by the ADMM-OtA-PINNs
achieve good agreements with the reference one and they cannot be visually distinguished. Moreover, as validated in Table \ref{tab:err_ex3}, we do not see much difference between the identified source by the ADMM-OtA-PINNs and the one by the ADMM-CG-FEM algorithm, which indicates that the ADMM-OtA-PINNs is comparable to traditional FEM-based numerical methods.

\begin{table}[h]
	\setlength{\abovecaptionskip}{0pt}
	\setlength{\belowcaptionskip}{3pt}
	\centering
	\caption{Numerical errors of ADMM-CG-FEM and  ADMM-OtA-PINNs for Example 3. Noise level $\delta=10\%$.}\label{tab:err_ex3}
	{\small\begin{tabular}{|c|c|c|c|}
			\hline
			Algorithm&ADMM-CG-FEM&ADMM-OtA-PINNs\\
			\hline
			$\frac{\|u-u_{true}\|}{\|u_{true}\|}$&0.175 &0.171 \\
			\hline
		\end{tabular}
	}
\end{table}


\section{Sparse Optimal Control of Parabolic Equations}\label{se:sparse_control}

In PDE-constrained optimal control problems, we usually can only put the controllers in some small regions instead of the whole domain under investigation. As a consequence, a natural question arises: how to determine the optimal locations and the intensities of the controllers? This concern inspires a class of optimal control problems where the controls are sparse (i.e., they are only non-zero in a small region of the domain); and the so-called sparse optimal control problems are obtained, see \cite{stadler2009elliptic}.  In this section,  we consider a parabolic sparse optimal control problem to validate the effectiveness and efficiency of the ADMM-PINNs algorithmic framework.

\subsection{Sparse Optimal Control Model}
A representative formulation of parabolic sparse optimal control problems can be formulated as:
 \begin{equation}\label{model_sparse}
	\underset{y,u\in L^2(Q) }{\min}~ \frac{1}{2}\|y-y_d\|^2_{L^2(Q)}+\frac{\alpha}{2}\|u\|^2_{L^2(Q)}+{\rho}\|u\|_{L^1(Q)}+I_{U_{ad}}(u),
 \end{equation}
where $y$ and $u$ satisfy the parabolic state equation:
\begin{flalign}\label{state_equation}
		\frac{\partial y}{\partial t}-\nu\Delta y+c_0y=u+f ~ \text{in}~\Omega\times(0,T), ~
		y=0~ \text{on}~ \partial\Omega\times(0,T),
		~y(0)=\varphi.
\end{flalign}
In (\ref{model_sparse})-(\ref{state_equation}), $\Omega$ is a bounded domain in $\mathbb{R}^d$ with $d \geq 1$ and $ \partial\Omega$ is its boundary; $Q=\Omega\times(0,T)$ with $0<T<+\infty$;  the desired state $y_d\in L^2(Q)$ and the functions $\phi\in L^2(\Omega), f\in L^2(Q)$ are given. The constants $\alpha>0$ and $\rho>0$ are regularization parameters. The coefficients $\nu>0$ and
$c_0\geq 0$ are assumed to be constant.
We denote by $I_{U_{ad}}(\cdot)$ the indicator function of the admissible set $U_{ad}$ defined by
\begin{equation*}\label{admissible}
	U_{ad}=\{u\in L^\infty(\Omega)| a\leq u(x,t)\leq b, ~\text{a.e.~in}~ Q\}\subset L^2(Q),
\end{equation*}
where $a,b \in L^2(\Omega)$ with $a < 0 < b$ almost everywhere.

Due to the presence of the nonsmooth $L^1$-regularization term, the structure of the optimal control of (\ref{model_sparse}) differs significantly from what one obtains for the usual smooth $L^2$-regularization. Particularly, the optimal control has small support as discussed in \cite{stadler2009elliptic,wachsmuth2011}. Because of this special structural property, optimal control problems with $L^1$-regularization capture important applications in various fields such as optimal actuator placement \cite{stadler2009elliptic} and impulse control \cite{ciaramella2016}. For example, as studied in \cite{stadler2009elliptic}, if one cannot or does not want to distribute control devices all over the control domain, but wants to place available devices in an optimal way, the solution of the problem (\ref{model_sparse}) can give some information about the optimal location of control devices.

\subsection{ADMM-PINNs for (\ref{model_sparse})}
Let $y(u)$ be the solution of the state equation (\ref{state_equation}) corresponding to $u$. By introducing $z\in L^2(Q)$ satisfying $u=z$, problem (\ref{model_sparse}) can be reformulated as
 {\small \begin{equation}\label{model_sparse_e}
		\begin{aligned}
		\underset{u, z\in L^2(Q) }{\min}~& \frac{1}{2}\|y(u)-y_d\|^2_{L^2(Q)}+\frac{\alpha}{2}\|u\|^2_{L^2(Q)}+{\rho}\|z\|_{L^1(Q)}+I_{U_{ad}}(z),\quad
		\text{s.t.}~&~u=z.
	\end{aligned}
 \end{equation}}
The augmented Lagrangian functional associated with (\ref{model_sparse_e}) is defined as
{\small
 \begin{equation}\label{L:sparse_e}
L_{\beta}(u, z;\lambda)=\frac{1}{2}\|y(u)-y_d\|^2_{L^2(Q)}+\frac{\alpha}{2}\|u\|^2_{L^2(Q)}+{\rho}\|z\|_{L^1(Q)}+I_{U_{ad}}(z)-(\lambda,u-z)_{L^2(Q)}+\frac{\beta}{2}\|u-z\|_{L^2(Q)}^2.
 \end{equation}
}
The ADMM algorithm for (\ref{model_sparse}) follows from \eqref{admm} with $L_{\beta}(u,z;\lambda)$ defined by \eqref{L:sparse_e}.

First, it is easy to show that the $z$-subproblem (\ref{admm_z}) reads as
 $$
z^{k+1}=\arg\min_{z\in L^2(Q)}I_{U_{ad}}(z)+\rho\|z\|_{L^1(Q)}-(\lambda^k,u^{k+1}-z)_{L^2(Q)}+\frac{\beta}{2}\|u^{k+1}-z\|_{L^2(Q)}^2,
 $$
which has a closed-form solution
$
z^{k+1}=\mathbb{P}_{U_{ad}}\left(\mathbb{S}_{\frac{\rho}{\beta}}\left(u^{k+1}-\frac{\lambda^k}{\beta}\right)\right).
$
Moreover, the $u$-subproblem (\ref{admm_u}) can be reformulated as
 \begin{equation}\label{u_sub_sparse}
\min_{y,u} \mathcal{J}_{SC}^k(y,u):=\frac{1}{2}\|y-y_d\|^2_{L^2(Q)}+\frac{\alpha}{2}\|u\|^2_{L^2(Q)}-(\lambda^k,u-z^k)_{L^2(Q)}+\frac{\beta}{2}\|u-z^k\|_{L^2(Q)}^2
 \end{equation}
subject to the state equation (\ref{state_equation}). The first-order optimality system of (\ref{u_sub_sparse}) is
\begin{equation*}
	\left\{
	\begin{aligned}
		&p +(\alpha+\beta) u-\lambda^k-\beta z^k=0,\\
	&\frac{\partial y}{\partial t}-\nu\Delta y+c_0y=u+f ~\text{in}~\Omega\times(0,T),\quad y=0~ \text{on}~ \partial\Omega\times(0,T),\quad y(0)=\varphi,\\
	&-\frac{\partial p}{\partial t}-\nu\Delta p+c_0p= y-y_d ~\text{in}~\Omega\times(0,T),\quad p=0~\text{on}~ \partial\Omega\times(0,T),\quad p(T)=0,
	\end{aligned}
	\right.
 \end{equation*}
where $u$ is the optimal solution of problem (\ref{u_sub_sparse}), $y$ and $p$ are the corresponding state and adjoint variable, respectively.

As discussed in Section \ref{se:solution_u}, the $u$-subproblem (\ref{admm_u}) can  be solved by applying Algorithm \ref{alg:pinn_d} or \ref{alg:pinn_oc}. Accordingly, two ADMM-PINNs algorithms can be obtained for solving problem (\ref{model_sparse}).  However, as mentioned in Section \ref{se:source_iden}, it was empirically observed that, when the ADMM-AtO-PINNs is applied, the corresponding neural networks are hard to train and the numerical results are not satisfactory. This is mainly because that the balance between the objective function and the PDE losses is sensitive and unstable for this case. To address this issue, much more effort is needed to find suitable loss weights, which is beyond the scope of the paper. Hence, we only focus on the application of the ADMM-OtA-PINNs to problem (\ref{model_sparse}).

\subsection{Numerical Results for (\ref{model_sparse})}
In this subsection, we report some numerical results to validate the effectiveness and efficiency of the ADMM-OtA-PINNs for solving the sparse optimal control problem (\ref{model_sparse}).

\medskip

\noindent\textbf{Example 4.} We consider an example of  problem (\ref{model_sparse}) with a known exact solution. In particular, we set $\Omega= (0, 1)^2, T = 1, \nu=1, c_0=0, a = -1,b = 2, \bar{y} =
5\sqrt{\rho}t\sin(3\pi x_1) \sin(\pi x_2), \bar{p} = 5\sqrt{\rho}(t - 1) \sin(\pi x_1)
\sin(\pi x_2)$, and
 $$
\bar{u}=\left\{
\begin{aligned}
&\max\{\frac{-\bar{p}+\rho}{\alpha},a\}~&&\text{in}~\{(x, t)\in \Omega\times (0,T): \bar{p}(x, t)>\rho\},\\
&\min\{\frac{-\bar{p}-\rho}{\alpha},b\}~&&\text{in}~\{(x, t)\in \Omega\times (0,T): \bar{p}(x, t)<-\rho\},\\
&0&&\text{otherwise}.
\end{aligned}
\right.
 $$
We further set $f=\frac{\partial \bar{y}}{\partial t}-\Delta \bar{y}-u$ and $y_d=\bar{y}-(-\frac{\partial \bar{p}}{\partial t}-\Delta\bar{p}).$
Then it can be shown that $\bar{u}$ is the optimal control of problem (\ref{model_sparse}) and $\bar{y}$ is the corresponding optimal state, the details can be referred to \cite{schiindele2017}

To solve the $u$-subproblem (\ref{admm_u}) in the PINNs framework, we approximate $y$ and $p$ with fully-connected feed-forward neural networks containing 3 hidden layers of 32 neurons each.  The hyperbolic tangent activation function is used in all the neural networks. To evaluate the loss, we uniformly sample $|\mathcal{T}_i|=4096$ residual points in the spatial-temporal domain $\Omega\times(0,T)$, and $|\mathcal{T}_{b_1}|=1024$ points in $\partial\Omega\times(0,T)$ and $|\mathcal{T}_{b_2}|=256$ points in $\Omega$ for the boundary and initial conditions.  Note that we reselect the residual points randomly at each ADMM iteration and they can be different from those in previous iterations. The weights are set as $w_y=w_p=1$, $w_i=1$ and $w_{b1}=w_{b2}=5$.

\begin{figure}[htpb]
	\caption{\small Numerical results for Example 4: sparsity in space.}\label{fig:control_sparse}
	\centering
		\subfloat[Exact control at $t=0.25$ ]{\includegraphics[width=0.32\textwidth]{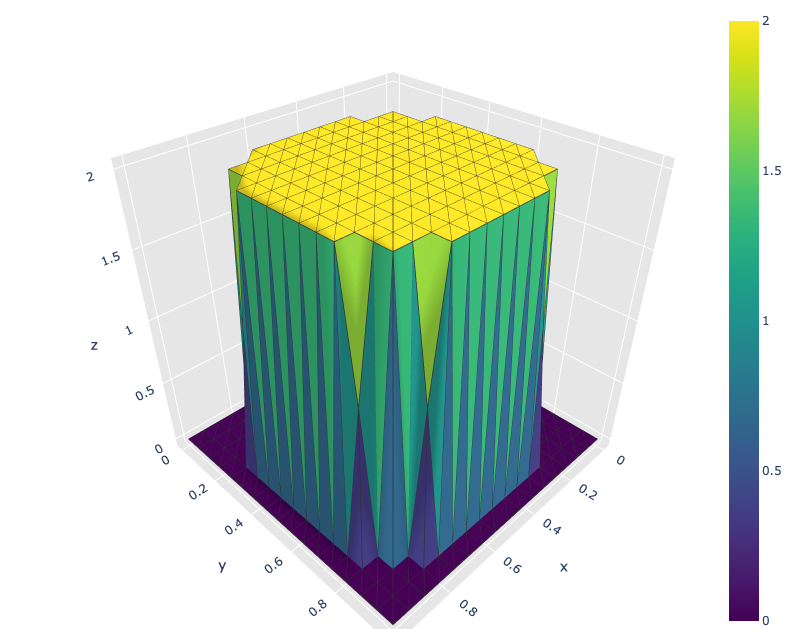}}
	\subfloat[Exact control at $t=0.5$]{\includegraphics[width=0.32\textwidth]{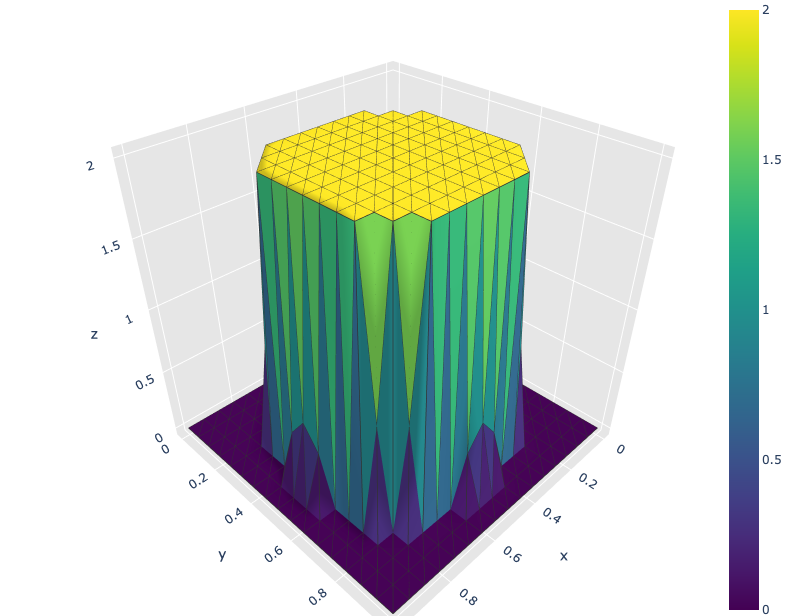}}
	\subfloat[Exact control at $t=0.75$]{\includegraphics[width=0.32\textwidth]{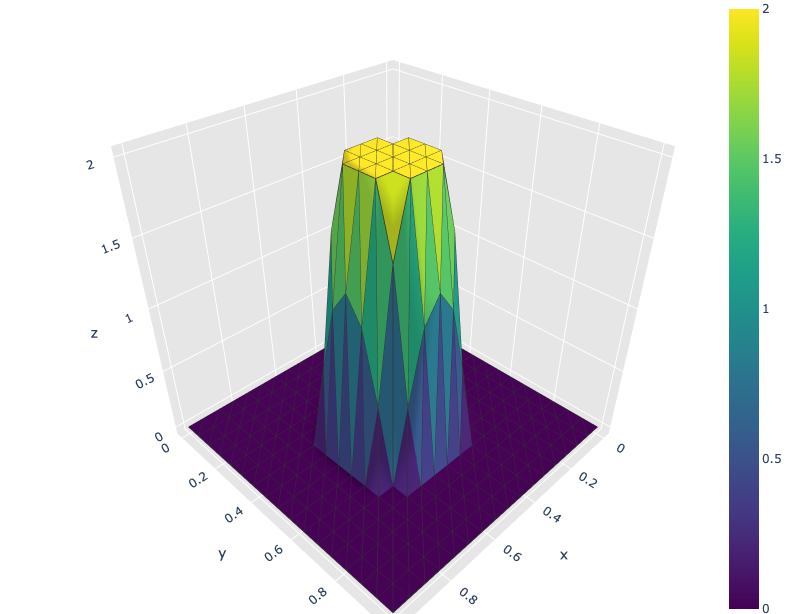}}\\
	\subfloat[Computed control at $t=0.25$]{\includegraphics[width=0.32\textwidth]{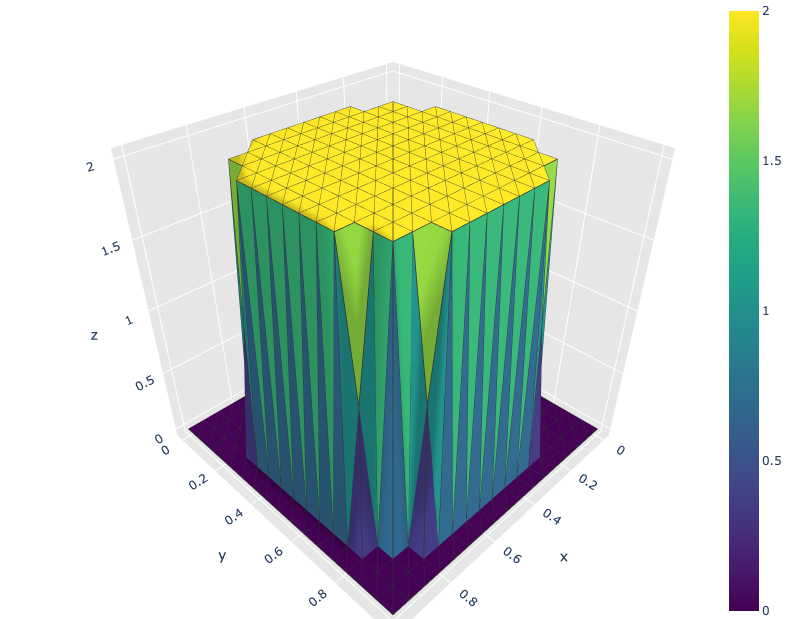}}
	\subfloat[Computed control  at $t=0.5$]{\includegraphics[width=0.32\textwidth]{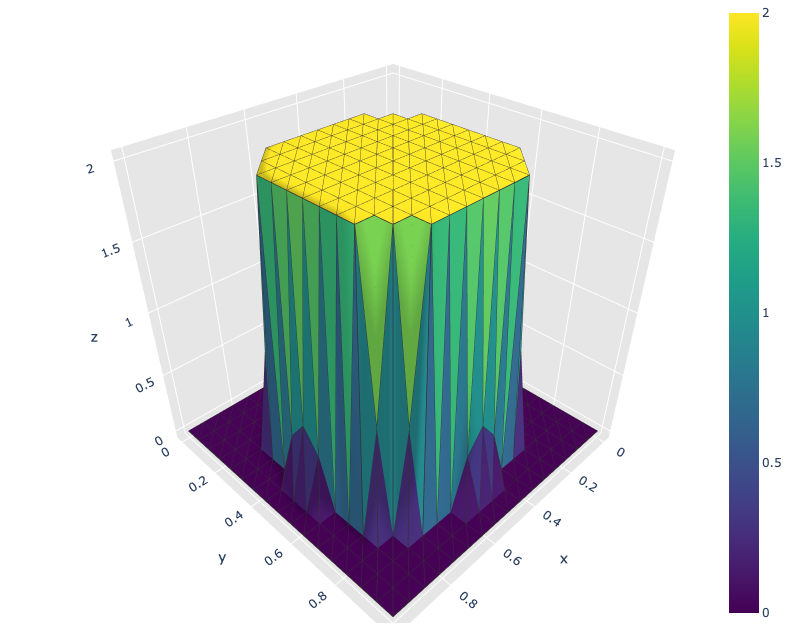}}
	\subfloat[Computed control at $t=0.75$]{\includegraphics[width=0.32\textwidth]{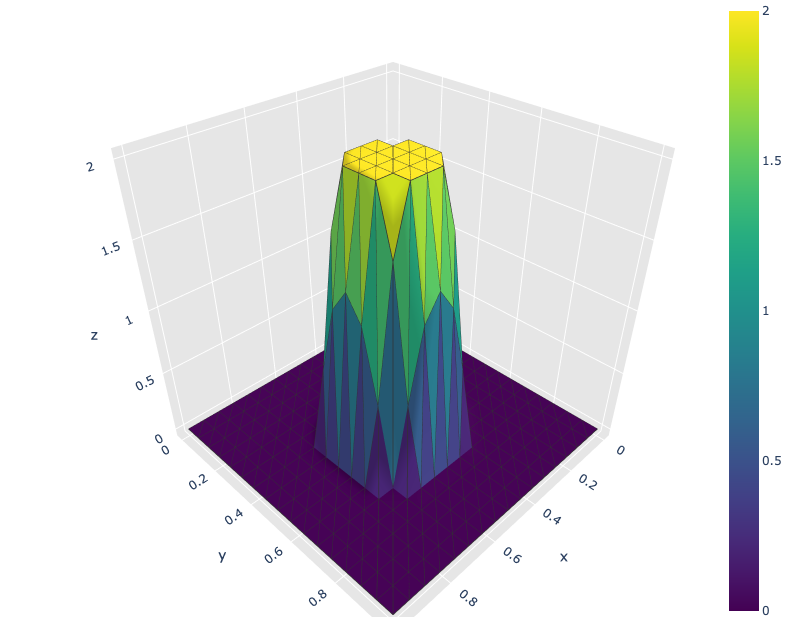}}
	\centering
\end{figure}

\begin{figure}[htpb]
	\caption{\small Numerical results for Example 4: sparsity in time.}\label{fig:time_sparse}
	\centering
	\includegraphics[width=\textwidth]{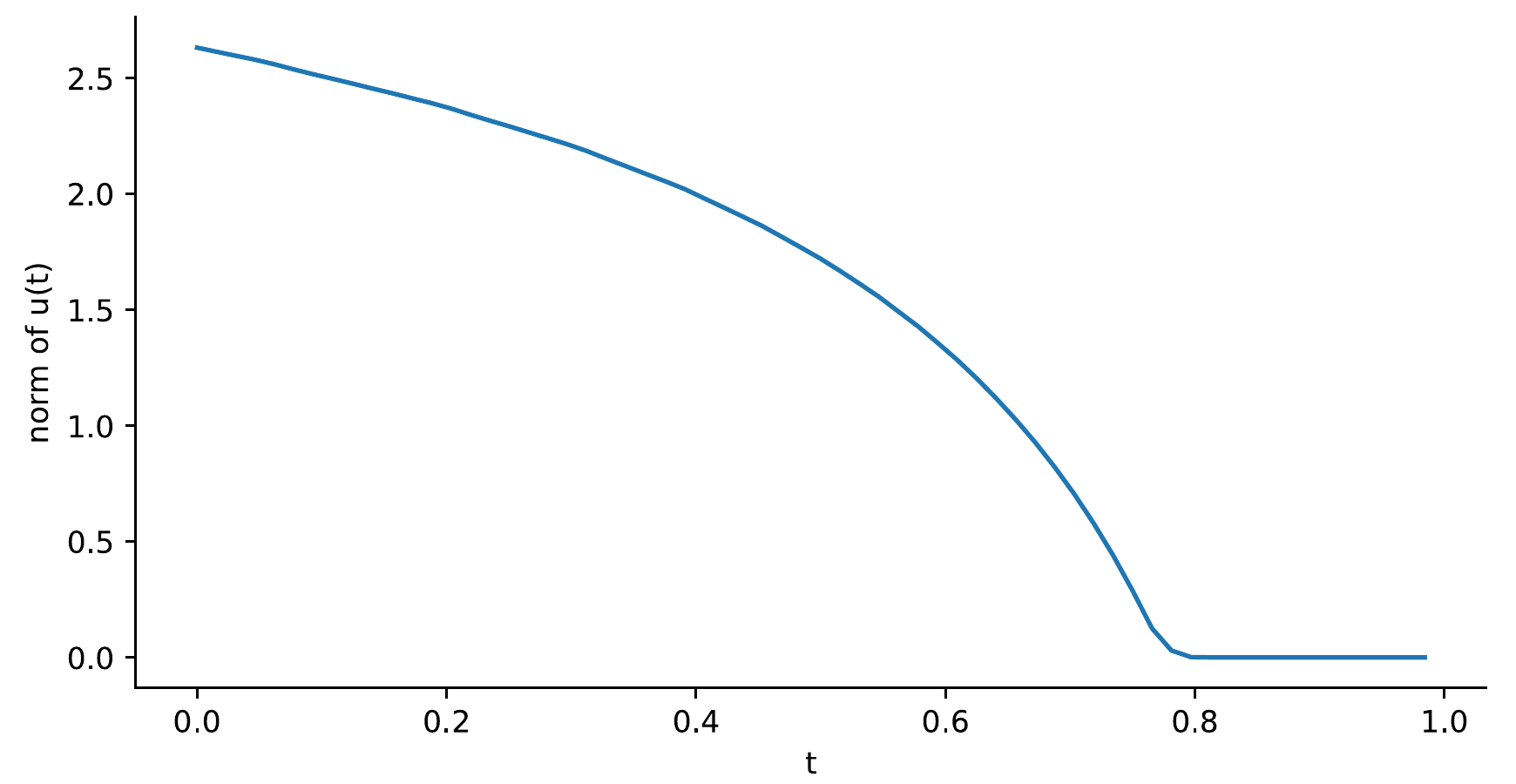}
\end{figure}

To train the neural networks,  we first use the Adam optimizer \cite{kingma2015} with learning rate $\eta=10^{-3}$ for 10000 iterations, and then switch to the L-BFGS \cite{nocedal1980} for 10 iterations.   In the implementation of the ADMM-OtA-PINNs, we execute 10 ADMM iterations with $\alpha=0.1, \rho=0.8,\beta=0.1$, $z^0=0$ and $\lambda^0=0$. The computed results are shown in Figures  \ref{fig:control_sparse} and \ref{fig:time_sparse}.

A comparison of the computed and exact controls at three different temporal snapshots ($t=0.25, 0.5$ and $0.75$)
is shown in Figure \ref{fig:control_sparse}. First, the sparsity of the control can be observed.  Moreover,  the computed control is very close to the exact one, and they even cannot be visually distinguished in Figure \ref{fig:control_sparse}.
Indeed, we have that the relative error
 $$
\frac{\|u^k(x, t)-\bar{u}(x, t)\|_{L^2(Q)}}{\|\bar{u}(x, t)\|_{L^2(Q)}}=1.45\%.
 $$

Furthermore, in our settings,  it is easy to see  $\bar{u}=0$ in $\{(x, t)\in \Omega\times (0,T): \bar{p}(x, t)<\rho\}$ and we can show that when $t>t^*=0.8211$, $u(x, t)=0$ a.e. in $\Omega$.  This conclusion is validated by the result presented in Figure \ref{fig:time_sparse}.
Hence, we conclude that a high-accuracy control can be pursued and the effectiveness of the proposed ADMM-PINNs algorithm is thus validated.

\section{Conclusions and Perspectives}\label{se:conclusion}

In this paper, we showcased a successful combination of contemporary deep learning techniques with classic operator splitting optimization techniques to design implementable and efficient algorithms for solving challenging PDE-constrained optimization problems. In particular, we suggested enhancing the well-known physics-informed neural networks (PINNs) with the classic alternating direction method of multipliers (ADMM), and thus significantly enlarged the applicable range of the PINNs to nonsmooth cases of PDE-constrained optimization problems. The proposed ADMM-PINNs algorithmic framework can be used to develop efficient algorithms for a class of nonsmooth PDE-constrained optimization problems, including inverse potential problems, source identification in elliptic equations, control constrained optimal control of the Burgers equation, and sparse optimal control of parabolic equations. The resulting algorithms differ from existing algorithms in that they do not require solving PDEs repeatedly and that they are mesh-free, easy to implement, and flexible to different PDEs. We expect the ADMM-PINNs algorithmic framework to pave an avenue of more application-tailored and sophisticated algorithms in the combination of traditional optimization techniques with emerging deep learning developments for solving more challenging problems arising in a broad array of science and engineering fields.

Our work leaves some important questions for future investigation, though solving these questions is beyond the scope of this paper. We list some of them below.

	$\bullet$ A natural extension of this work is to consider other even more challenging PDEs such as the Cahn–Hilliard equation \cite{hintermuller2012}, the Navier Stokes equation \cite{hintermuller2006} and the Boussinesq equation \cite{song2022}.

	$\bullet$ The validated efficiency of the ADMM-PINNs algorithms clearly justifies the necessity to study the underlying theoretical issues such as the convergence analysis and the error estimate. Note that the convergence of ADMM for convex problems has been extensively studied in the classic literature (see e.g., \cite{glowinski1984}). The generic model \eqref{pco_model} under our discussion, however, can be nonconvex and there is no unified framework for the convergence analysis of ADMM for nonconvex cases. On top of this, theoretical analysis of PINNs {\it per se} is drawing extensive attention but there are still no mature results for the generic model \eqref{pco_model}. Hence, convergence analysis for the ADMM-PINNs may require case-by-case studies for different problems.
	
	$\bullet$ The reliability of the ADMM-PINNs in different contexts crucially depends on the quantification of the total uncertainty associated with noisy data, neural network hyperparameters, and optimization and sampling errors. In this regard, some uncertainty quantification (UQ) methods proposed in the literature, e.g., \cite{psaros2023,yang2021,zhang2019,zou2023}, for PINNs will be useful. Moreover, the NeuralUQ library developed in \cite{zou2024} will help facilitate the numerical implementations of these UQ methods.

	$\bullet$ In our numerical experiments, Algorithm \ref{alg:pinn_d} or \ref{alg:pinn_oc} for solving the $u$-subproblem (\ref{u_sub}) in the ADMM-PINNs algorithm framework was terminated empirically. It has been shown in \cite{GSYY2022} that it is not necessary to pursue a too high-precision solution for (\ref{u_sub}), especially when the iterations are still far from the solution point. Also, solutions with too low precision for (\ref{u_sub}) may not be sufficient to guarantee overall convergence. Therefore, it is important to determine an adaptive inaccuracy criterion to implement the inner PINNs methods to solve (\ref{u_sub}).
	
   $\bullet$ In Sections \ref{se:source_iden} and \ref{se:sparse_control}, we observed that it is not easy to train the neural networks for the resulting ADMM-AtO-PINNs. The main reason is that the balance between the objective function and the PDE losses is sensitive and unstable. Similar results have also been observed in \cite{haoBilevel2022}. To further improve the numerical performance, it seems that some adaptive strategies should be developed for discerning the loss weights. Research on this topic seems still to be in its infancy.

\section*{Acknowledgments}
The authors would like to thank Professor Enrique Zuazua for helpful discussions and supports to the completion of this paper.  The authors are grateful to two anonymous referees for their very valuable comments which have helped improve the paper substantially.

\bibliographystyle{siamplain}

\end{document}